\documentclass[11pt,twoside]{amsart} % {report}{article}{book}{amsart} [twoside, landscape ,draft]
\usepackage{amssymb}
\catcode`@=11

\let\latex@newcommand\newcommand
\let\latex@renewcommand\renewcommand

\def\if@undefined#1{\ifx#1\un@@@defined@@@}
\def\newcommand#1{\if@undefined{#1}
  \let\next@=\latex@newcommand \else
  \let\next@=\latex@renewcommand \fi
  \next@{#1}}

\catcode`@=12
\theoremstyle{definition}

    \newtheorem{thm}{Theorem}[section]
    \newtheorem{prop}[thm]{Proposition}
    \newtheorem{lemma}[thm]{Lemma}
    \newtheorem{defn}[thm]{Definition}
    \newtheorem{cor}[thm]{Corollary}

    \newtheorem*{thm*}{Theorem}
    \newtheorem*{prop*}{Proposition}
    \newtheorem*{lemma*}{Lemma}
    \newtheorem*{defn*}{Definition}
    \newtheorem*{cor*}{Corollary}
    \newtheorem{rems}[thm]{Remark}
    \newtheorem*{rems*}{Remark}
    \newtheorem*{proof*}{Proof}
    \newtheorem*{not*}{Notation}

\renewenvironment{proof}{{\it Proof.}}{$\Box$\par}

\newcommand{\ndef}{\newcommand}
\def\rndef{\renewcommand}

\ndef{\myaddress}[1]{\begin{center} \it\tiny #1 \end{center}}

\ndef{\clA}{{\mathcal A}} \ndef{\rmA}{{\mathrm A}} \ndef{\mbA}{{\mathbb A}} \ndef{\bfA}{{\mathbf A}} \ndef{\euA}{{\EuScript A}} \ndef{\frA}{{\mathfrak A}}
\ndef{\clB}{{\mathcal B}} \ndef{\rmB}{{\mathrm B}} \ndef{\mbB}{{\mathbb B}} \ndef{\bfB}{{\mathbf B}} \ndef{\euB}{{\EuScript B}} \ndef{\frB}{{\mathfrak B}}
\ndef{\clC}{{\mathcal C}} \ndef{\rmC}{{\mathrm C}} \ndef{\mbC}{{\mathbb C}} \ndef{\bfC}{{\mathbf C}} \ndef{\euC}{{\EuScript C}} \ndef{\frC}{{\mathfrak C}}
\ndef{\clD}{{\mathcal D}} \ndef{\rmD}{{\mathrm D}} \ndef{\mbD}{{\mathbb D}} \ndef{\bfD}{{\mathbf D}} \ndef{\euD}{{\EuScript D}} \ndef{\frD}{{\mathfrak D}}
\ndef{\clE}{{\mathcal E}} \ndef{\rmE}{{\mathrm E}} \ndef{\mbE}{{\mathbb E}} \ndef{\bfE}{{\mathbf E}} \ndef{\euE}{{\EuScript E}} \ndef{\frE}{{\mathfrak E}}
\ndef{\clF}{{\mathcal F}} \ndef{\rmF}{{\mathrm F}} \ndef{\mbF}{{\mathbb F}} \ndef{\bfF}{{\mathbf F}} \ndef{\euF}{{\EuScript F}} \ndef{\frF}{{\mathfrak F}}
\ndef{\clG}{{\mathcal G}} \ndef{\rmG}{{\mathrm G}} \ndef{\mbG}{{\mathbb G}} \ndef{\bfG}{{\mathbf G}} \ndef{\euG}{{\EuScript G}} \ndef{\frG}{{\mathfrak G}}
\ndef{\clH}{{\mathcal H}} \ndef{\rmH}{{\mathrm H}} \ndef{\mbH}{{\mathbb H}} \ndef{\bfH}{{\mathbf H}} \ndef{\euH}{{\EuScript H}} \ndef{\frH}{{\mathfrak H}}
\ndef{\clI}{{\mathcal I}} \ndef{\rmI}{{\mathrm I}} \ndef{\mbI}{{\mathbb I}} \ndef{\bfI}{{\mathbf I}} \ndef{\euI}{{\EuScript I}} \ndef{\frI}{{\mathfrak I}}
\ndef{\clJ}{{\mathcal J}} \ndef{\rmJ}{{\mathrm J}} \ndef{\mbJ}{{\mathbb J}} \ndef{\bfJ}{{\mathbf J}} \ndef{\euJ}{{\EuScript J}} \ndef{\frJ}{{\mathfrak J}}
\ndef{\clK}{{\mathcal K}} \ndef{\rmK}{{\mathrm K}} \ndef{\mbK}{{\mathbb K}} \ndef{\bfK}{{\mathbf K}} \ndef{\euK}{{\EuScript K}} \ndef{\frK}{{\mathfrak K}}
\ndef{\clL}{{\mathcal L}} \ndef{\rmL}{{\mathrm L}} \ndef{\mbL}{{\mathbb L}} \ndef{\bfL}{{\mathbf L}} \ndef{\euL}{{\EuScript L}} \ndef{\frL}{{\mathfrak L}}
\ndef{\clM}{{\mathcal M}} \ndef{\rmM}{{\mathrm M}} \ndef{\mbM}{{\mathbb M}} \ndef{\bfM}{{\mathbf M}} \ndef{\euM}{{\EuScript M}} \ndef{\frM}{{\mathfrak M}}
\ndef{\clN}{{\mathcal N}} \ndef{\rmN}{{\mathrm N}} \ndef{\mbN}{{\mathbb N}} \ndef{\bfN}{{\mathbf N}} \ndef{\euN}{{\EuScript N}} \ndef{\frN}{{\mathfrak N}}
\ndef{\clO}{{\mathcal O}} \ndef{\rmO}{{\mathrm O}} \ndef{\mbO}{{\mathbb O}} \ndef{\bfO}{{\mathbf O}} \ndef{\euO}{{\EuScript O}} \ndef{\frO}{{\mathfrak O}}
\ndef{\clP}{{\mathcal P}} \ndef{\rmP}{{\mathrm P}} \ndef{\mbP}{{\mathbb P}} \ndef{\bfP}{{\mathbf P}} \ndef{\euP}{{\EuScript P}} \ndef{\frP}{{\mathfrak P}}
\ndef{\clQ}{{\mathcal Q}} \ndef{\rmQ}{{\mathrm Q}} \ndef{\mbQ}{{\mathbb Q}} \ndef{\bfQ}{{\mathbf Q}} \ndef{\euQ}{{\EuScript Q}} \ndef{\frQ}{{\mathfrak Q}}
\ndef{\clR}{{\mathcal R}} \ndef{\rmR}{{\mathrm R}} \ndef{\mbR}{{\mathbb R}} \ndef{\bfR}{{\mathbf R}} \ndef{\euR}{{\EuScript R}} \ndef{\frR}{{\mathfrak R}}
\ndef{\clS}{{\mathcal S}} \ndef{\rmS}{{\mathrm S}} \ndef{\mbS}{{\mathbb S}} \ndef{\bfS}{{\mathbf S}} \ndef{\euS}{{\EuScript S}} \ndef{\frS}{{\mathfrak S}}
\ndef{\clT}{{\mathcal T}} \ndef{\rmT}{{\mathrm T}} \ndef{\mbT}{{\mathbb T}} \ndef{\bfT}{{\mathbf T}} \ndef{\euT}{{\EuScript T}} \ndef{\frT}{{\mathfrak T}}
\ndef{\clU}{{\mathcal U}} \ndef{\rmU}{{\mathrm U}} \ndef{\mbU}{{\mathbb U}} \ndef{\bfU}{{\mathbf U}} \ndef{\euU}{{\EuScript U}} \ndef{\frU}{{\mathfrak U}}
\ndef{\clV}{{\mathcal V}} \ndef{\rmV}{{\mathrm V}} \ndef{\mbV}{{\mathbb V}} \ndef{\bfV}{{\mathbf V}} \ndef{\euV}{{\EuScript V}} \ndef{\frV}{{\mathfrak V}}
\ndef{\clW}{{\mathcal W}} \ndef{\rmW}{{\mathrm W}} \ndef{\mbW}{{\mathbb W}} \ndef{\bfW}{{\mathbf W}} \ndef{\euW}{{\EuScript W}} \ndef{\frW}{{\mathfrak W}}
\ndef{\clX}{{\mathcal X}} \ndef{\rmX}{{\mathrm X}} \ndef{\mbX}{{\mathbb X}} \ndef{\bfX}{{\mathbf X}} \ndef{\euX}{{\EuScript X}} \ndef{\frX}{{\mathfrak X}}
\ndef{\clY}{{\mathcal Y}} \ndef{\rmY}{{\mathrm Y}} \ndef{\mbY}{{\mathbb Y}} \ndef{\bfY}{{\mathbf Y}} \ndef{\euY}{{\EuScript Y}} \ndef{\frY}{{\mathfrak Y}}
\ndef{\clZ}{{\mathcal Z}} \ndef{\rmZ}{{\mathrm Z}} \ndef{\mbZ}{{\mathbb Z}} \ndef{\bfZ}{{\mathbf Z}} \ndef{\euZ}{{\EuScript Z}} \ndef{\frZ}{{\mathfrak Z}}

\ndef{\tA}{{\widetilde A}} \ndef{\tcA}{{\widetilde\clA}} \ndef{\ttcA}{\widetilde{\tcA}} \ndef{\sfA}{{\textsf A}} \ndef{\ttA}{\widetilde{\tA}} \ndef{\dzA}{{A^sharp}}
\ndef{\tB}{{\widetilde B}} \ndef{\tcB}{{\widetilde\clB}} \ndef{\ttcB}{\widetilde{\tcB}} \ndef{\sfB}{{\textsf B}} \ndef{\ttB}{\widetilde{\tB}} \ndef{\dzB}{{B^sharp}}
\ndef{\tC}{{\widetilde C}} \ndef{\tcC}{{\widetilde\clC}} \ndef{\ttcC}{\widetilde{\tcC}} \ndef{\sfC}{{\textsf C}} \ndef{\ttC}{\widetilde{\tC}} \ndef{\dzC}{{C^sharp}}
\ndef{\tD}{{\widetilde D}} \ndef{\tcD}{{\widetilde\clD}} \ndef{\ttcD}{\widetilde{\tcD}} \ndef{\sfD}{{\textsf D}} \ndef{\ttD}{\widetilde{\tD}} \ndef{\dzD}{{D^sharp}}
\ndef{\tE}{{\widetilde E}} \ndef{\tcE}{{\widetilde\clE}} \ndef{\ttcE}{\widetilde{\tcE}} \ndef{\sfE}{{\textsf E}} \ndef{\ttE}{\widetilde{\tE}} \ndef{\dzE}{{E^sharp}}
\ndef{\tF}{{\widetilde F}} \ndef{\tcF}{{\widetilde\clF}} \ndef{\ttcF}{\widetilde{\tcF}} \ndef{\sfF}{{\textsf F}} \ndef{\ttF}{\widetilde{\tF}} \ndef{\dzF}{{F^sharp}}
\ndef{\tG}{{\widetilde G}} \ndef{\tcG}{{\widetilde\clG}} \ndef{\ttcG}{\widetilde{\tcG}} \ndef{\sfG}{{\textsf G}} \ndef{\ttG}{\widetilde{\tG}} \ndef{\dzG}{{G^sharp}}
\ndef{\tH}{{\widetilde H}} \ndef{\tcH}{{\widetilde\clH}} \ndef{\ttcH}{\widetilde{\tcH}} \ndef{\sfH}{{\textsf H}} \ndef{\ttH}{\widetilde{\tH}} \ndef{\dzH}{{H^sharp}}
\ndef{\tI}{{\widetilde I}} \ndef{\tcI}{{\widetilde\clI}} \ndef{\ttcI}{\widetilde{\tcI}} \ndef{\sfI}{{\textsf I}} \ndef{\ttI}{\widetilde{\tI}} \ndef{\dzI}{{I^sharp}}
\ndef{\tJ}{{\widetilde J}} \ndef{\tcJ}{{\widetilde\clJ}} \ndef{\ttcJ}{\widetilde{\tcJ}} \ndef{\sfJ}{{\textsf J}} \ndef{\ttJ}{\widetilde{\tJ}} \ndef{\dzJ}{{J^sharp}}
\ndef{\tK}{{\widetilde K}} \ndef{\tcK}{{\widetilde\clK}} \ndef{\ttcK}{\widetilde{\tcK}} \ndef{\sfK}{{\textsf K}} \ndef{\ttK}{\widetilde{\tK}} \ndef{\dzK}{{K^sharp}}
\ndef{\tL}{{\widetilde L}} \ndef{\tcL}{{\widetilde\clL}} \ndef{\ttcL}{\widetilde{\tcL}} \ndef{\sfL}{{\textsf L}} \ndef{\ttL}{\widetilde{\tL}} \ndef{\dzL}{{L^sharp}}
\ndef{\tM}{{\widetilde M}} \ndef{\tcM}{{\widetilde\clM}} \ndef{\ttcM}{\widetilde{\tcM}} \ndef{\sfM}{{\textsf M}} \ndef{\ttM}{\widetilde{\tM}} \ndef{\dzM}{{M^sharp}}
\ndef{\tN}{{\widetilde N}} \ndef{\tcN}{{\widetilde\clN}} \ndef{\ttcN}{\widetilde{\tcN}} \ndef{\sfN}{{\textsf N}} \ndef{\ttN}{\widetilde{\tN}} \ndef{\dzN}{{N^sharp}}
\ndef{\tO}{{\widetilde O}} \ndef{\tcO}{{\widetilde\clO}} \ndef{\ttcO}{\widetilde{\tcO}} \ndef{\sfO}{{\textsf O}} \ndef{\ttO}{\widetilde{\tO}} \ndef{\dzO}{{O^sharp}}
\ndef{\tP}{{\widetilde P}} \ndef{\tcP}{{\widetilde\clP}} \ndef{\ttcP}{\widetilde{\tcP}} \ndef{\sfP}{{\textsf P}} \ndef{\ttP}{\widetilde{\tP}} \ndef{\dzP}{{P^sharp}}
\ndef{\tQ}{{\widetilde Q}} \ndef{\tcQ}{{\widetilde\clQ}} \ndef{\ttcQ}{\widetilde{\tcQ}} \ndef{\sfQ}{{\textsf Q}} \ndef{\ttQ}{\widetilde{\tQ}} \ndef{\dzQ}{{Q^sharp}}
\ndef{\tR}{{\widetilde R}} \ndef{\tcR}{{\widetilde\clR}} \ndef{\ttcR}{\widetilde{\tcR}} \ndef{\sfR}{{\textsf R}} \ndef{\ttR}{\widetilde{\tR}} \ndef{\dzR}{{R^sharp}}
\ndef{\tS}{{\widetilde S}} \ndef{\tcS}{{\widetilde\clS}} \ndef{\ttcS}{\widetilde{\tcS}} \ndef{\sfS}{{\textsf S}} \ndef{\ttS}{\widetilde{\tS}} \ndef{\dzS}{{S^sharp}}
\ndef{\tT}{{\widetilde T}} \ndef{\tcT}{{\widetilde\clT}} \ndef{\ttcT}{\widetilde{\tcT}} \ndef{\sfT}{{\textsf T}} \ndef{\ttT}{\widetilde{\tT}} \ndef{\dzT}{{T^sharp}}
\ndef{\tU}{{\widetilde U}} \ndef{\tcU}{{\widetilde\clU}} \ndef{\ttcU}{\widetilde{\tcU}} \ndef{\sfU}{{\textsf U}} \ndef{\ttU}{\widetilde{\tU}} \ndef{\dzU}{{U^sharp}}
\ndef{\tV}{{\widetilde V}} \ndef{\tcV}{{\widetilde\clV}} \ndef{\ttcV}{\widetilde{\tcV}} \ndef{\sfV}{{\textsf V}} \ndef{\ttV}{\widetilde{\tV}} \ndef{\dzV}{{V^sharp}}
\ndef{\tW}{{\widetilde W}} \ndef{\tcW}{{\widetilde\clW}} \ndef{\ttcW}{\widetilde{\tcW}} \ndef{\sfW}{{\textsf W}} \ndef{\ttW}{\widetilde{\tW}} \ndef{\dzW}{{W^sharp}}
\ndef{\tX}{{\widetilde X}} \ndef{\tcX}{{\widetilde\clX}} \ndef{\ttcX}{\widetilde{\tcX}} \ndef{\sfX}{{\textsf X}} \ndef{\ttX}{\widetilde{\tX}} \ndef{\dzX}{{X^sharp}}
\ndef{\tY}{{\widetilde Y}} \ndef{\tcY}{{\widetilde\clY}} \ndef{\ttcY}{\widetilde{\tcY}} \ndef{\sfY}{{\textsf Y}} \ndef{\ttY}{\widetilde{\tY}} \ndef{\dzY}{{Y^sharp}}
\ndef{\tZ}{{\widetilde Z}} \ndef{\tcZ}{{\widetilde\clZ}} \ndef{\ttcZ}{\widetilde{\tcZ}} \ndef{\sfZ}{{\textsf Z}} \ndef{\ttZ}{\widetilde{\tZ}} \ndef{\dzZ}{{Z^sharp}}

\ndef{\bfc}{{\bf c}}

\let\geq\geqslant
\let\leq\leqslant

\ndef{\lims}[1]{\lim\limits_{#1}}
\ndef{\sums}[1]{\sum\limits_{#1}}
\ndef{\ints}[1]{\int\limits_{#1}}
\ndef{\sups}[1]{\sup\limits_{#1}}
\ndef{\liminfty}[1]{\lims{#1\to\infty}}
\ndef{\suminf}[1]{\sums{#1=1}^\infty}

\ndef{\limo}[1]{\omega\mbox{-}\!\!\!\lims{#1\to\infty}}
\ndef{\limL}[1]{\rmL\mbox{-}\!\!\!\lims{#1\to\infty}}
\ndef{\limLOne}[1]{\clL_1\mbox{-}\!\!\!\lims{#1}}
\ndef{\tildelimo}[1]{\tilde\omega\mbox{-}\!\!\!\lims{#1\to\infty}}

\ndef{\normE}[1]{\norm{#1}_E}               % the norm of #1 in a symmetric Banach space E
\ndef{\Aut}{\operatorname{Aut}}
\ndef{\Ch}{\operatorname{ch}}        % Chern character
\ndef{\End}{\operatorname{End}}
\ndef{\Hom}{\operatorname{Hom}}
\ndef{\Ker}{\operatorname{Ker}}
\ndef{\Log}{\operatorname{Log}}
\ndef{\OP}{\operatorname{OP}}
\ndef{\Op}{\operatorname{Op}}
\ndef{\Symb}{\operatorname{Symb}}
\ndef{\Tr}{\operatorname{Tr}}
\ndef{\Wres}{\operatorname{Wres}}
\ndef{\cl}{\operatorname{cl}}
\ndef{\com}{\operatorname{com}}
\ndef{\const}{\operatorname{const}}
\ndef{\conv}{\operatorname{conv}}
\rndef{\det}{\operatorname{det}}
\ndef{\detFK}{\operatorname{det_{FK}}}
\ndef{\diag}{\operatorname{diag}}
\ndef{\dist}{\operatorname{dist}}
\ndef{\dom}{\operatorname{dom}}
\ndef{\ec}{\operatorname{ec}}        % essential codimension
\ndef{\id}{1}
\ndef{\ind}{\operatorname{ind}}
\ndef{\mydeg}{\operatorname{deg}}
\ndef{\op}{\operatorname{op}}
\ndef{\rank}{\operatorname{rank}}
\ndef{\res}{\operatorname{res}}      % residue
\ndef{\rng}{\operatorname{ran}}      % range
\ndef{\sflow}{\operatorname{sf}}     % spectral flow
\ndef{\sign}{\operatorname{sign}}
\ndef{\sing}{\operatorname{sing}}
\ndef{\supp}{\operatorname{supp}}
\ndef{\tr}{\operatorname{tr}}
\ndef{\vol}{\operatorname{vol}}      % volume or volume form
\ndef{\wn}{\operatorname{wn}}        % winding number
\ndef{\wres}{\operatorname{wres}}    % Wodzicki residue
\rndef{\Im}{\operatorname{Im}}
\rndef{\Re}{\operatorname{Re}}

\ndef{\rslv}[1]{R_z(#1)}
\ndef{\HH}{H}
\ndef{\tHH}{\tilde \HH}
\ndef{\VV}{V}
\ndef{\Rz}{R_z}
\ndef{\tRz}{\tR_z}
\ndef{\psif}[1]{#1^{[1]}} % {\psi_{#1}}
\ndef{\bndl}{\xi}                         % vector bundle
\ndef{\bndlA}{\eta}                      % vector bundle
\ndef{\GlueMap}{\varphi}                        % glue map of a bundle
\ndef{\ChartMap}{h}                             % chart diffeomorphism map of a manifold

\ndef{\hilb}{\clH}                     % hilbert space
   \ndef{\hilbasargument}{(\hilb)} %{(\hilb)}
\ndef{\LpH}[1]{\clL^{#1}}                      % the set of ...
\ndef{\TrCl}{\LpH{1}}                          % the set of TRACE-CLASS operators
\ndef{\TrClH}{\LpH{1}\hilbasargument}          % the set of TRACE-CLASS operators on \hilb
\ndef{\clBH}{\clB\hilbasargument}              % the set of BOUNDED linear operators on hilbert space
\ndef{\clCH}{\clC\hilbasargument}              % the set of CLOSED DENSELY-DEFINED linear operators on hilbert space
\ndef{\clKH}{\clK\hilbasargument}              % the set of COMPACT operators
\ndef{\clFH}{\clF\hilbasargument}              % the set of BOUNDED FREDHOLM operators
\ndef{\clUH}{\clU\hilbasargument}              % the set of UNITARIES on hilbert space
\ndef{\clCFH}{{\clC\clF}\hilbasargument}       % the set of CLOSED DENSELY-DEFINED FREDHOLM OPERATORS on hilbert space
\ndef{\saBH}{\clB\hilbasargument_{sa}}         % the set of S.-A. BOUNDED operators on hilbert space
\ndef{\saCH}{\clC\hilbasargument_{sa}}         % the set of CLOSED DENSELY-DEFINED S.-A. operators on hilbert space
\ndef{\saFH}{\clF\hilbasargument_{sa}}         % the set of BOUNDED FREDHOLM s.-a. operators
\ndef{\saKH}{\clK\hilbasargument_{sa}}         % the set of COMPACT S.-A. operators
\ndef{\saCFH}{\clC\clF\hilbasargument_{sa}}    % the set of CLOSED DENSELY-DEFINED S.-A. FREDHOLM operators on hilbert space
\ndef{\clUFH}{\clU\clF\hilbasargument}         % the set of UNITARIES such that U+I is FREDHOLM
\ndef{\Uinj}{\clU_{inj}\hilbasargument}        % the set of UNITARIES such that U-I is injective
\ndef{\UFinj}{\clU\clF_{inj}\hilbasargument}   % the set of UNITARIES such that U-I is INJECTIVE and U+I is FREDHOLM

\ndef{\LpN}[1]{\clL^{#1}(\clN,\tau)}       % noncommutative \mathcal L_p space
\ndef{\rLpN}[1]{L^{#1}(\clN,\tau)}       % noncommutative L_p space
\ndef{\clAND}{(\clA,\clN,\clD)}          % spectral triple (A,N,D)
\ndef{\clBA}{{\clB(\clA)}}
\ndef{\clKN}{{\clK(\clN,\tau)}}          % \tau-compact operators
\ndef{\clKtN}{{\clK_\tau(\tilde\clN)}}   % \tau-compact (maybe unbounded) operators
\ndef{\clFN}{{\clF_\tau(\clN)}}          % \tau-Fredholm operators
\ndef{\clPN}{\clP(\clN)}                 % projections of N
\ndef{\clQN}{\clQ_\tau(\clN)}            % Calkin algebra N/K
\ndef{\infPN}{{\clP_\tau^\infty(\clN)}}  % infinite projections of N
\ndef{\clOF}[2]{\clF_{#1\mbox{-}#2}(\clN,\tau)}         % relatively Fredholm operators
\ndef{\oind}[2]{{\rm \tau\mbox{-}ind}_{#1\mbox{-}#2}}   % relative index
\ndef{\tind}{\ind_\tau}                  % semifinite index
\ndef{\DInd}{\ind_{\clD,\tau}}           % semifinite index
\ndef{\BF}{$\tau$-Fredholm }
\ndef{\affl}{\eta}
\ndef{\vNa}{von Neumann algebra}
\ndef{\nsf}{faithful semifinite normal }  % normal semifinite faithful
\ndef{\taubrs}[1]{\tau\brackets{#1}}

\ndef{\domd}{\bigcap\limits_{n\ge 0} \dom\;\delta^n}          % domain of \delta^n's
\ndef{\DiffOP}{{\rm \clD}}
\ndef{\ADA}{\clA \cup [\clD,\clA]}
\ndef{\DixIdeal}[1]{\LpH{#1,\infty}}               % Dixmier ideal
\ndef{\dixideal}{\ell^{1,\infty}}                  % Dixmier ideal
\ndef{\WDixIdeal}{\LpH{1,\mathrm w}}               % weak Dixmier ideal
\ndef{\DixIdealPos}[1]{\DixIdeal{#1}_+}          % positive part of Dixmier ideal
\ndef{\DixIdealN}[1]{\LpN{#1,\infty}}            % semifinite Dixmier ideal
\ndef{\DixIdealNPar}[2]{\clL_{#1,\infty}^{#2}(\clN,\tau)}    % semifinite Dixmier ideal
\ndef{\DixIdealNPos}[1]{\LpN{#1,\infty}_+}       % positive part of semifinite Dixmier ideal
\ndef{\TrD}{\Tr_\omega}                                       % Dixmier trace
\ndef{\tauD}{{\tau_\omega}}                                   % semifinite Dixmier trace
\ndef{\ILog}{\frac 1{\log(1+t)}}
\ndef{\ILogN}{\frac 1{\log(1+N)}}
\ndef{\DixNorm}[1]{\norm{#1}_{(1,\infty)}}                        % Dixmier norm
\ndef{\DixInt}[1]{\ints 0^t \mu_s(#1)\,ds}
\ndef{\DixIntL}[1]{\ints 0^{\lambda_{1/t}(#1)}\mu_s(#1)\,ds}
    \ndef{\SmallIdeal}{\clL_{1, \mathrm w}}
    \ndef{\DixIntII}[1]{\ints 0^t \mu_s(#1)ds}
    \ndef{\DixIntf}[1]{f_t(#1)}
    \ndef{\DixIntg}[1]{g_t(#1)}

\ndef{\lpi}{\LpN{1,(\pi)}}
\ndef{\HaarMeasBohrs}{\nu}            % Haar measure of Bohr's compact
\ndef{\BrownsMeas}{\mu}               % Brown's measure
\ndef{\BohrCont}[1]{\tilde{#1}}        % continuation on Bohr's compact
\ndef{\APMean}{{M}}                   % mean value of a.p. function
\ndef{\CDSS}{{\clA_B}}                % Coburn-Douglas-Schaeffer-Singer's factor
\ndef{\matr}{{\rm Mat}}
\ndef{\seque}[1]{\ensuremath{\{#1_j\}_{j=1}^\infty}}    % sequence of numbers  a_1, a_2, ...
\ndef{\sequen}[2]{\ensuremath{\{#1_#2\}_{#2=1}^\infty}}    % sequence of numbers  a_1, a_2, ...
\ndef{\Seque}[1]{\ensuremath{\left(#1_0,#1_1,#1_2,\dots\right)}}    % sequence of numbers  a_1, a_2, ...
\ndef{\Cesaro}{H}                           % the Cesaro operator (on sequences)
\ndef{\CesaroRPlus}{M}                      % the Cesaro operator on positive semiaxis
\ndef{\Dilation}{D}                         % the dilation operator (on sequences)
\ndef{\Shift}{T}                            % the shift operator (on sequences)
\ndef{\norm}[1]{\left\Vert#1\right\Vert}          % norm of #1
\ndef{\TrNorm}[1]{\norm{#1}_1}              % trace norm of #1
\ndef{\HSNorm}[1]{\norm{#1}_2}              % Hilbert-Schmidt norm of #1
\ndef{\InftyNorm}[1]{\norm{#1}_\infty}      % uniform norm of #1
\ndef{\abs}[1]{\left\lvert#1\right\rvert}   % absolute value of #1
\ndef{\set}[1]{\left\{#1\right\}}           % set of ...
\ndef{\brackets}[1]{\left(#1\right)}
\ndef{\brs}[1]{\brackets{#1}}
\ndef{\scalprod}[2]{\la #1,#2\ra}
\ndef{\precprec}{\prec\!\!\!\prec}
\ndef{\qeq}{\stackrel?=}
\ndef{\spectrum}[1]{\sigma_{#1}} %{\mathrm{Spec}(#1)}
\ndef{\numrange}[1]{\mathrm{W}(#1)}
\rndef{\emptyset}{\varnothing}
\ndef{\sss}[1]{\subsubsection{}\label{#1}}
\rndef{\phi}{\varphi}
\ndef{\OpenUnitDisk}{D}
\ndef{\RHS}{RHS}
\ndef{\LHS}{LHS} %right and left hand side
\ndef{\ttt}{\Leftrightarrow}
\ndef{\then}{\Rightarrow}
\ndef{\tto}{\longrightarrow}
\ndef{\nno}{\nonumber\\}
\ndef{\newn}[1]{\index{#1} \emph{#1}}                                % new notion
\ndef{\la}{\langle}
\ndef{\ra}{\rangle}
\ndef{\dbar}{{\;\bar{\phantom{o}} \!\!\!\! d}}
\ndef{\stl}[1]{\stackrel{\vbox to 0pt{\vss\hbox{$\scriptstyle #1$}}}}
\ndef{\mathcomment}[1]{{\scriptstyle\text{(#1)}}\qquad}        % for comments at the ends of lines with math formulas
\ndef{\details}[1]{\smallskip\begin{center} {\bf Here:} #1\end{center}\medskip}

     \ndef{\npartial}{\slash\!\!\!\partial}
     \ndef{\Heis}{\operatorname{Heis}}
     \ndef{\Solv}{\operatorname{Solv}}
     \ndef{\Spin}{\operatorname{Spin}}
     \ndef{\SO}{\operatorname{SO}}
     \ndef{\Index}{\operatorname{index}}
             \ndef{\coker}{{\mbox coker}}
             \ndef{\p}{\partial}
             \ndef{\dd}{|\clD|}
             \ndef{\n}{\parallel}
     \ndef{\gf}[2]{\genfrac{}{}{0pt}{}{#1}{#2}}
     \ndef{\ta}{\widetilde{\alpha}}
     \ndef{\tb}{\widetilde{\beta}}
     \ndef{\txi}{\widetilde{\xi}}
     \ndef{\tk}{\widetilde{K}}
     \ndef{\CGh}{\widetilde{\CG}}
     \ndef{\boe}{{\bf e}}\ndef{\bt}{{\bf t}}
     \ndef{\vth}{\vartheta}
     \ndef{\db}{\overline{\partial}}
     \ndef{\hV}{\hat{V}}
     \ndef{\cag}{{\clA^\Gamma}}
     \ndef{\sind}{\sigma{\rm -ind}}

\begin{document}
\theoremstyle{definition}

\def\clA{{\mathcal A}}    \def\rmA{{\mathrm A}}    \def\mbA{{\mathbb A}}
\def\clB{{\mathcal B}}    \def\rmB{{\mathrm B}}    \def\mbB{{\mathbb B}}
\def\clC{{\mathcal C}}    \def\rmC{{\mathrm C}}    \def\mbC{{\mathbb C}}
\def\clD{{\mathcal D}}    \def\rmD{{\mathrm D}}    \def\mbD{{\mathbb D}}
\def\clE{{\mathcal E}}    \def\rmE{{\mathrm E}}    \def\mbE{{\mathbb E}}
\def\clF{{\mathcal F}}    \def\rmF{{\mathrm F}}    \def\mbF{{\mathbb F}}
\def\clG{{\mathcal G}}    \def\rmG{{\mathrm G}}    \def\mbG{{\mathbb G}}
\def\clH{{\mathcal H}}    \def\rmH{{\mathrm H}}    \def\mbH{{\mathbb H}}
\def\clI{{\mathcal I}}    \def\rmI{{\mathrm I}}    \def\mbI{{\mathbb I}}
\def\clJ{{\mathcal J}}    \def\rmJ{{\mathrm J}}    \def\mbJ{{\mathbb J}}
\def\clK{{\mathcal K}}    \def\rmK{{\mathrm K}}    \def\mbK{{\mathbb K}}
\def\clL{{\mathcal L}}    \def\rmL{{\mathrm L}}    \def\mbL{{\mathbb L}}
\def\clM{{\mathcal M}}    \def\rmM{{\mathrm M}}    \def\mbM{{\mathbb M}}
\def\clN{{\mathcal N}}    \def\rmN{{\mathrm N}}    \def\mbN{{\mathbb N}}
\def\clO{{\mathcal O}}    \def\rmO{{\mathrm O}}    \def\mbO{{\mathbb O}}
\def\clP{{\mathcal P}}    \def\rmP{{\mathrm P}}    \def\mbP{{\mathbb P}}
\def\clQ{{\mathcal Q}}    \def\rmQ{{\mathrm Q}}    \def\mbQ{{\mathbb Q}}
\def\clR{{\mathcal R}}    \def\rmR{{\mathrm R}}    \def\mbR{{\mathbb R}}
\def\clS{{\mathcal S}}    \def\rmS{{\mathrm S}}    \def\mbS{{\mathbb S}}
\def\clT{{\mathcal T}}    \def\rmT{{\mathrm T}}    \def\mbT{{\mathbb T}}
\def\clU{{\mathcal U}}    \def\rmU{{\mathrm U}}    \def\mbU{{\mathbb U}}
\def\clV{{\mathcal V}}    \def\rmV{{\mathrm V}}    \def\mbV{{\mathbb V}}
\def\clW{{\mathcal W}}    \def\rmW{{\mathrm W}}    \def\mbW{{\mathbb W}}
\def\clX{{\mathcal X}}    \def\rmX{{\mathrm X}}    \def\mbX{{\mathbb X}}
\def\clY{{\mathcal Y}}    \def\rmY{{\mathrm Y}}    \def\mbY{{\mathbb Y}}
\def\clZ{{\mathcal Z}}    \def\rmZ{{\mathrm Z}}    \def\mbZ{{\mathbb Z}}
%%%%%%%%%%%%%%%%%%%NEW MACROS%%%%%%%%%%%%%%%%%%%%%%%
\def\clNt{\widetilde{\clN}}
\def\nt{(\clN,\tau)}
\def\ent{E\nt}
\def\clMt{\widetilde{\clM}}
\def\LL{{\mathbf L}}
\def\R{{\mathbf R}}
\def\frB{{\mathfrak B}}
%%%%%%%%%%%%%%%%%%%%%%%%%%%%%%%%%

\def\lims#1{\lim\limits_{#1}}
\def\sums#1{\sum\limits_{#1}}
\def\ints#1{\int\limits_{#1}}
\def\sups#1{\sup\limits_{#1}}
\def\liminfty#1{\lim\limits_{#1\to\infty}}

\def\limo#1{\omega\mbox{-}\!\!\!\lim\limits_{#1\to\infty}} %\mathsurround 0pt
\def\limL#1{\rmL\mbox{-}\!\!\!\lim\limits_{#1\to\infty}} %\mathsurround 0pt
\def\tildelimo#1{\tilde\omega\mbox{-}\!\!\!\lim\limits_{#1\to\infty}} %\mathsurround 0pt

%        NG notation
%%%%%%%%%%%%%%%%%%%%%%%%%%%%%%%%%%%%%%%%%%%%%%%%%%%%%%%%%%%%%%%%%%%%%%%%%%%%%%%%%%%%%%%%%%%%%%%%%%%%%%%%%%%%%
\def\domd{\bigcap\limits_{n\ge 0} \dom\;\delta^n}          % domain of \delta^n's
\def\DiffOP{{\rm \clD}}
\def\ADA{\clA \cup [\clD,\clA]}
\def\DixIdeal#1{\clL^{(#1,\infty)}(\hilb)}               % Dixmier ideal
\def\DixIdealPos#1{\clL_+^{(#1,\infty)}(\hilb)}          % positive part of Dixmier ideal
\def\DixIdealN#1{\clL^{(#1,\infty)}(\clN,\tau)}            % semifinite Dixmier ideal
\def\DixIdealNPar#1#2{\clL^{({#1},\infty)}_{#2}(\clN,\tau)}    % semifinite Dixmier ideal
\def\DixIdealNPos#1{\clL^{(#1,\infty)}_+(\clN,\tau)}       % positive part of semifinite Dixmier ideal
\def\tauD{{\tau_\omega}}                                   % semifinite Dixmier trace
\def\ILog{\frac 1{\log(1+t)}}
\def\ILogN{\frac 1{\log(1+N)}}
\def\DixNorm#1{||#1||_{(1,\infty)}}                        % Dixmier norm
\def\DixInt#1{\ints 0^t \mu_s(#1)\,ds}
\def\DixIntL#1{\ints 0^{\l_{1/t}(#1)}\mu_s(#1)\,ds}

%%%%%%%%%%%%%%%%%%%%%%%%%%%%%%%%%%%%%%%%%%%%%%%%%%%%%%%%%%%%%%%%%%%%%%%%%%%%%%%%%%%%%%%%%%%%%%%%%%%%%%%%%%%%%
%        General mathematical notation
%%%%%%%%%%%%%%%%%%%%%%%%%%%%%%%%%%%%%%%%%%%%%%%%%%%%%%%%%%%%%%%%%%%%%%%%%%%%%%%%%%%%%%%%%%%%%%%%%%%%%%%%%%%%%
\def\seque#1{\ifmmode \{#1_j\}_{j=1}^\infty \else $\{#1_j\}_{j=1}^\infty$ \fi}    % sequence of numbers  a_1, a_2, ...
\def\sequen#1#2{\ifmmode \{#1_#2\}_{#2=1}^\infty \else $\{#1_#2\}_{#2=1}^\infty$ \fi}    % sequence of numbers  a_1, a_2, ...
\def\Seque#1{\ifmmode \left(#1_1,#1_2,#1_3,\dots\right) \else $\left(#1_1,#1_2,#1_3,\dots\right)$ \fi}    % sequence of numbers  a_1, a_2, ...
\def\Cesaro{H}                     % the Cesaro operator (on sequences)
\def\CesaroRPlus{M}                     % the Cesaro operator on positive semiaxis
\def\Dilation{D}                   % the dilation operator (on sequences)
\def\Shift{T}                                 % the shift operator (on sequences)
\def\then{\Rightarrow}
\def\newn#1{\index{#1} \emph{#1}}
\def\la{\langle} \def\ra{\rangle}

\def\title#1{\begin{center}\bf\large #1\end{center}\vskip 0.3 in}
\def\author#1{\begin{center}#1\end{center}}

\titlepage
\begin{center}
  \title{THE DIXMIER TRACE AND ASYMPTOTICS OF ZETA FUNCTIONS}
\vspace{.2 in}
{\bf by}\\
\vspace{.2 in}
{\bf Alan L. Carey}$^{1,3}$\\Mathematical Sciences Institute\\
Australian National University\\
Canberra ACT 0200\\AUSTRALIA\\
\vspace{.2 in}
{\bf Adam Rennie}$^2$\\
Institute for Mathematical Sciences\\University of Copenhagen\\
Universitetsparken 5, DK-2100, Copenhagen, DENMARK\\
\vspace{.2 in}
{\bf Aleksandr Sedaev}$^{1,4}$\\
Department of Mathematics,\\ Voronezh State
University of Architecture and Construction,\\
 20-letiya Oktyabrya
84, Voronezh 394006, RUSSIA\\
\vspace{.2 in}
{\bf Fyodor Sukochev}$^1$\\School of Informatics and Engineering\\
Flinders University\\Bedford Park 5042\\AUSTRALIA\\
\vspace{0.25 in}
AMS classification nos.
 Primary: 19K56, 46L80: secondary: 58B30, 46L87.\\
Keywords:
Spectral triple, Dixmier trace, zeta function.

\vspace{0.25 in}
Supported by the Clay Mathematics Institute$^3$
and by grants from ARC$^1$ (Australia), SNF (Denmark)$^2$ and RFBR (5-01-00629)(Russia)$^4$.

%\end{center}
{\bf Abstract}

\end{center}
We obtain general theorems which enable the calculation of the
Dixmier trace in terms of the asymptotics of the zeta function
and of the trace of the heat semigroup.
We prove our results in a general semi-finite von Neumann
algebra. We find for $p>1$ that the
asymptotics of the zeta function determines
an ideal strictly larger than ${\mathcal L}^{p,\infty}$ on which
the Dixmier trace may be defined. We also establish stronger versions of
other results on Dixmier traces and zeta functions.
%\end{abstract}
%\maketitle
%\end{center}
%\tableofcontents
\newpage
\section{Introduction}

\subsection{Background}
The key role of the Dixmier trace in noncommutative geometry was
discovered by Connes around 1990, \cite{Co4}. Since then, it has
become a cornerstone of noncommutative geometry. Notably, the
Dixmier trace is used to define dimension, integration and has
been used in physical applications, along with heat kernel type
expansions, to define `spectral actions' for noncommutative field
theories, \cite{chc,Co6}. The Dixmier trace (or more precisely
Dixmier {\em traces}) are a family of non-normal traces on
the bounded operators on a separable Hilbert space $\mathcal H$
measuring the logarithmic divergence
of the trace of a compact operator. There is an ideal of compact
operators denoted  ${\mathcal
  L}^{(1,\infty)}({\mathcal H})$ consisting precisely of those
operators with finite
Dixmier trace. (This and the related ideals  ${\mathcal
  L}^{(p,\infty)}({\mathcal H})$, $p\geq 1$,
are defined in detail in  Section 2 {\it cf} also \cite{Co4}.)
Following \cite{Co4}
connections between Dixmier traces, zeta functions and heat kernel
asymptotics were systematically studied in \cite{CPS2}. Motivated by
these results, and questions arising in connection with physical
applications, we substantially extend the understanding of these matters
 in this article.

Briefly, for an important special case, we show that for a
positive compact operator $T$,
the existence of the limit
$\lim_{r\to\infty}\frac{1}{r}\mbox{Trace}(T^{p+\frac{1}{r}})
$
implies that the operator $T$ lies in an ideal ${\mathcal
Z}_p$. The ideal
${\mathcal Z}_1$ is $\clL^{(1,\infty)}({\mathcal H})$, while for $p>1$ ${\mathcal Z}_p$ is strictly larger than
 $\clL^{(p,\infty)}({\mathcal H})$.
(It is  in fact precisely
what is termed, in \cite[Section 1.d]{LT2},
 the $p$-convexification of $\clL^{(1,\infty)}({\mathcal H})$.)
We then show that if
$\lim_{r\to\infty}\frac{1}{r}\mbox{Trace}(T^{p+\frac{1}{r}})$ exists it equals
$p\mbox{Trace}_\omega(T^p) \label{intro-two}$
for any state $\omega$ generating a Dixmier trace, Trace$_\omega$.
Thus we show  that the asymptotics
 of the zeta function
singles out the class of compact operators
which have a finite Dixmier trace.

In fact the analogues of these statements
are true for compact operators $T$ in a semifinite
von Neumann algebra ${\mathcal N}$ with faithful, normal, semifinite
trace $\tau$ for which there are
corresponding ideals ${\mathcal Z}_p({\mathcal N})$ and
$\clL^{(p,\infty)}({\mathcal N},\tau)$. Readers unfamiliar with ideal
theory in such general algebras
may restrict attention to the standard case of bounded operators on
an infinite dimensional  separable Hilbert space with its usual trace
(denoted by `Trace' here).
Our reason for striving for generality stems from
the emergence recently of applications of the
semifinite von Neumann theory \cite{BCPRSW, BeF, CM,  CPRS2, CPRS3, pr}.

Our results follow primarily
from (strengthened versions of) deep facts from  \cite{CPS2}
and recent advances in the study of singular traces, some of which
seem not to be well known.
We also work in this paper with general
Marcinkiewicz spaces and general `Dixmier traces'
as these spaces are already known to arise
in the study of pseudodifferential operators \cite{N}.

Before giving a more precise account of our results, let us set out
the motivations coming from noncommutative physics and geometry.
In \cite{GGISV} it
was shown that the Moyal `plane' of dimension $2N$ defines a
$(2N,\infty)$-summable spectral triple. In order to prove this, the
authors used a variant of Cwikel's inequality, and to compute Dixmier
traces, they employed the zeta function methods of
\cite{CPS2}. Numerous other noncommutative spaces which are
$(p,\infty)$-summable have been studied, \cite{cp1,cp2,chc,CoL,DLSSV,pr,prs}, some with physical
applications or relevance.

Examining these examples shows
that except for very special and/or simple examples, eg
\cite{cp1,cp2,pr,prs}, the determination of Dixmier summability of an operator
relies on one of two methods: Weyl's theorem, or Cwikel type
inequalities. In particular for operators arising from `noncommutative
action principles' (that is when we minimise functionals on noncommutative
algebras), no (classical) geometric context need exist, and so
Weyl's theorem is of no use.

The theorems presented here offer alternative techniques for proving
Dixmier summability results, and computing Dixmier traces. This is
likely to be relevant for (very) noncommutative examples and
physically inspired examples. It is also likely that via zeta function
regularisation of determinants, our techniques could provide criteria
for one-loop renormalizability of noncommutative field
theories.
%Finally, our results can be used to study
%summability criteria for nonunital spectral triples, at least when a
%suitable local structure is assumed, \cite{r}.

\subsection{Summary of the main results}
We need some notation in order to present the results.
We remark that in a semifinite von Neumann algebra ${\mathcal N}$ with faithful
normal semifinite trace $\tau$ the  $\tau-$compact
operators are generated by projections $P$ with $\tau(P)<\infty$. Suppose
that $T$ is a $\tau-$compact positive operator in ${\mathcal N}$.
(If one has a semifinite
spectral triple determined by an unbounded self adjoint operator
$D$ then one should think of $T$ as $|D|^{-1}$ or
$(1+D^2)^{-1/2}$.)
For a given $\tau$  let $\tau_\omega$ denote a Dixmier trace corresponding to
an element $\omega\in \ell_\infty^*(\mathbb N)$
or $\ell_\infty^*(\mathbb R_+)$. We remark that
$\omega$ must satisfy some invariance properties
which we will explain in detail in Section 3.
%Our intention is that the inclusion of
% these Sections will serve to clarify
%the different notions of singular trace which appear in the
%literature.
By the zeta function of $T$ we mean $\zeta(s)=\tau(T^s)$.

Consider the following hypothesis:\\
%(i) Under the assumption that $e^{-tT^{-2}}$ is trace class for all $t>0$
%suppose that
%$$\lim_{\lambda\to\infty}\lambda^{-1}\tau(Ae^{-\lambda^{-2/p}T^{-2}})$$
%$$\omega-\lim_{\lambda\to\infty}\lambda^{-1}\tau(Ae^{-\lambda^{-2}T^{-2}})$$
%exists where $\omega$ is a dilation and Cesaro mean invariant functional
%on $L^\infty({\mathbb R}^+)$.
(*) Under the assumption that $\tau(T^{s})$
exists for all $s>p$ suppose that
$
\lim_{r\to\infty} \frac{1}{r}\zeta(p+\frac{1}{r})$
exists.\\
It is then natural to ask, in view of \cite{CPS2,Co4}, the following question:\\A. If hypothesis (*) holds then does it follow that $T\in \clL ^{(p,\infty)}$?\\We prove that the answer to Question A is yes if $p=1$ and no if $p>1$.
This leads to a second question:\\
B. For $p>1$ what constraint does hypothesis (*) place on the singular values
of $T$?

We remark that in contrast to the situation with the classical
Schatten ideals it is not true that if $T\in \clL ^{(1,\infty)}$
then $T^{1/p}\in \clL ^{(p,\infty)}$. In fact there is a strictly
smaller ideal inside $\clL ^{(1,\infty)}$ characterized by this
property. We prove correspondingly that there is an ideal
${\mathcal Z}_p$ strictly larger than $\clL ^{(p,\infty)}$ with
the property that if $T^{1/p}\in {\mathcal Z}_p$ then $T\in \clL
^{(1,\infty)}$. We also prove that if hypothesis (*) holds then
$T\in {\mathcal Z}_p$.

This leads to the further question:\\
C. If hypothesis (*) holds how does the limit
relate to the Dixmier trace of $T^p$?\\
In fact  we show that for a certain class of Dixmier traces $\tau_\omega$ 
$$\lim_{r\to\infty} \frac{1}{r}\zeta(p+\frac{1}{r})=p\lim_{t\to\infty}\frac1{\log(1+t)}\int_0^t\mu_s(T)^pds:=p\tau_\omega(T^p).$$

Our methods then lead us to prove some stronger versions of several results in
\cite{CPS2}. In that paper we were forced to consider a subset of
the set of all Dixmier traces determined by requiring invariance under
a certain transformation group. In the
new approach of this article we can relax many of these invariance conditions.

Then, in view of \cite[p 563]{Co4} and the relationship of the zeta
function to the heat kernel, it is natural to ask what hypothesis
(*) implies concerning the small time asymptotics of the trace of
the heat semigroup. (We note that hypothesis (*) implies that the
heat semigroup $e^{-tT^{-2}}$, defined using the functional
calculus, is trace class for all $t>0$.) This matter is resolved
in Theorem \ref{mainprop}. Let
$F(\lambda)=\lambda^{-1}\tau(e^{-\lambda^{-2}T^{-2}})$, then under
hypothesis (*) for $p=1$ this function is bounded on $(0,\infty)$
and Theorem \ref{mainprop} says that for certain
$\omega\in L_\infty((0,\infty))^*$, $\omega(F)$ is
a multiple of the Dixmier trace $\tau_\omega(T)$.

%$$\lim_{\lambda\to\infty}\lambda^{-1}\tau(Ae^{-\lambda^{-2/p}T^{-2}})$$
%$$\omega-\lim_{\lambda\to\infty}\lambda^{-1}\tau(Ae^{-\lambda^{-2}T^{-2}})$$
%exists where $\omega$ is a dilation and Cesaro mean invariant functional
%on $L^\infty({\mathbb R}^+)$.

Conversely we know
that if $\lambda^{-1}\tau(e^{-\lambda^{-2}T^{-2}})$
has an asymptotic expansion in $\lambda$ as $\lambda\to \infty$
then the leading term in this expansion precisely determines the 
first singularity of $\tau(T^s)$ as $\Re(s)$ decreases. In this case,
using the results described above, we find that
$T\in{\mathcal Z}_1$ 
and the residue of the zeta function is equal to the Dixmier trace of $T$.

Finally, in Section 6, we revisit a question
raised in \cite{CGS}. Namely,  for $T$ in some general ideal $\mathcal I$
 (in the $\tau-$compact operators),
which admits a Dixmier trace $\tau_\omega$,
what are the minimal conditions on
an algebra $\mathcal A$ such that the functional $a\to
\tau_\omega(aT)$ on $\mathcal A$ is actually a trace? This question is
important in the manifold reconstruction theorem of \cite{Co4}. We
find that the methods of this paper enable us to substantially
generalize \cite{CGS} (who answer the question only for
 ${\mathcal I}=\clL^{(p,\infty)}$). We find that,
for the same minimal conditions as in \cite{CGS}, there is a very
large class of Marcinkiewicz
ideals $\mathcal I$ including ${\mathcal I}={\mathcal Z}_p$
for which $a\to
\tau_\omega(aT)$ is a trace.

We  give in Section 2 a summary of the theory of singular 
traces and a careful discussion of ideals of
compact operators needed in this paper. We follow this
in Section 3 with some
details on the  construction of Dixmier traces. 
The main results are proved in Section 4, for the
zeta function, and Section 5 for the heat operator. We finish 
with our generalization of \cite{CGS}.

\noindent{\bf Acknowledgement} The fourth named author
 thanks Bruno Iochum for asking
the question that led to Corollary \ref{bruno}
and the Universit\'e de Cergy-Pontoise for hosting his visit. We also thank
Harald Grosse and Victor Gayral for advice and Evgenii Semenov for
explaining to us some subtle facts about geometry of Marcinkiewicz
spaces. The first named author thanks the Erwin Schr\"odinger
Institute for its assistance with this research and the 
Clay Mathematics Institute for financial support.

%%%%%%%%%%%%%%%%%%%%end insertion%%%%%%%%%%%%%%%%%%
\let\le\leq
\let\ge\geq
\def\x{x}
\def\y{y}
\def\z{z}

%%%%%%%%%%%%%%%%%%%%%%%%%%%%%Kalton Sukochev I%%%%%%%%%%%%%%%%%%%%%%%%

{

\setlength{\belowdisplayshortskip}{1\belowdisplayskip}
\newcommand{\no}{\hspace*{\indentation}}
\renewcommand{\arraystretch}{1.75}
\setlength{\arrayrulewidth}{0.5\arrayrulewidth}

\newcommand{\tab}[1]{\hspace*{#1\tablength}}
\newcommand{\nm}[1]{\mbox{\ensuremath{\| #1 \|}}}
\newcommand{\implies}{\ensuremath{\Rightarrow}}
\newcommand{\fa}{\ensuremath{\ \, \forall \,}}
\renewcommand{\iff}{\ensuremath{\; \Longleftrightarrow \;}}
\renewcommand{\sp}[1]{\ensuremath{\mathrm{sp}(#1)}}
\newcommand{\gap}[1]{\ensuremath{\, #1 \,}}
\newcommand{\inset}[2]{\ensuremath{\{ #1 \, | \, #2 \} }}
\newcommand{\inprod}[2]{\ensuremath{\langle #1 , #2 \rangle}}
\newcommand{\anti}[2]{\ensuremath{ \{ #1 , #2 \} }}
\newcommand{\CC}{\ensuremath{\mathbb{C}}}
\newcommand{\RR}{\ensuremath{\mathbb{R}}}
\newcommand{\TT}{\ensuremath{\mathbb{T}}}
\newcommand{\SB}{\ensuremath{\mathbb{S}}}
\newcommand{\ZZ}{\ensuremath{\mathbb{Z}}}
\newcommand{\HH}{\ensuremath{\mathbb{H}}}
\newcommand{\NN}{\ensuremath{\mathbb{N}}}
\newcommand{\iny}{\ensuremath{\infty}}
\newcommand{\bd}[1]{ \ensuremath{r_{\NN}(\ensuremath{#1})}}
\newcommand{\bp}{\ensuremath{p}}
\newcommand{\nmt}{\widetilde {\mathcal {M}}}
\newcommand{\nmm}{\mathcal {M}}
\newcommand{\emt}{E(\mathcal {M},\tau )}

\newcommand{\noline}{\vspace*{1\parskip}}
\newcommand{\preskip}{\vspace*{1\belowdisplayskip}}
\newcommand{\display}[1]{$$#1$$}
\newcommand{\smdisplay}[1]{\\[4pt] \no \centerline{#1} \\[4pt]}
\newcommand{\text}[1]{\mbox{\normalfont #1}}
\newcommand{\mod}[1]{\ensuremath{
        \text{ \hspace*{-3pt} mod \hspace*{-3pt} } #1}}

\section{Preliminaries: spaces and functionals}
\subsection{Function spaces} \label{PrelimOnSpAndF}
The theory of singular traces on operator ideals rests on some classical
analysis which we now review for completeness.

\noindent Consider a Banach space $(E, \Vert \cdot \Vert_E)$ of
real valued Lebesgue measurable functions on the interval $J=[0,\infty)$
or else on
$J=\mbN$. Let $x^*$ denote the non-increasing, right-continuous
rearrangement of $|x|$ given by \display{x^*(t) = \inf \inset{s
\geq 0}{\lambda(\{ |x|> s\}) \leq t} , \ t > 0,} where $\lambda$
denotes Lebesgue measure. Then $E$ will be called rearrangement
invariant (or r.i.) if

\noindent (i).\quad $E$ is an ideal lattice, that is if
$y\in E$, and $x$ is any measurable function on $J$ with $0\leq
|x|\leq |y|$, then $x\in E$ and $\Vert x\Vert_E\leq \Vert
y\Vert_E$;\\
\noindent (ii).\quad if $y\in E$ and
 if $x$ is any
measurable function on $J$ with $x^*=y^*$, then $x\in E$ and
$\Vert x\Vert_E= \Vert y\Vert_E$.

\noindent In the case $J=\mbN$,
it is convenient to identify $x^*$ with the rearrangement of the
sequence $|x|=\{|x_n|\}_{n=1}^\infty$ in descending order.
(The theory is in the monographs \cite{KPS}, \cite{LT}, \cite{LT2}.)
A r.i. space $E$ is said to be a fully symmetric Banach space if
it has the additional property that if $y\in E$ and $L_1 +
L_\infty(J)\ni x\prec\prec y$, then $x\in E$ and $\Vert
x\Vert_E\leq \Vert y\Vert_E$. Here, $x\prec\prec y$ denotes
submajorization in the sense of Hardy-Littlewood-P\'olya:
$$\int_0^tx^*(s)ds\leq  \int_0^ty^*(s)ds,\quad \fa t>0.
$$
All these spaces $E$ satisfy $L_{1}\cap L_{\infty }\left( J\right)
\subseteq E\subseteq \left( L_{1}+L_{\infty }\right) \left( J
\right) $, with continuous embeddings. In this paper, we consider
only fully symmetric Banach spaces $E$, which satisfy in addition
$E\subseteq L_{\infty }(J)$ (a non-commutative extension of the
theory of such spaces placed in the setting of a semifinite von
Neumann algebra ${\mathcal N}$ corresponds to ideals in ${\mathcal
N}$ equipped with unitarily invariant norm
\cite{GK,DDP,CS,Suk,CPS2}).

Recall (see  \cite{KPS}) that for an arbitrary rearrangement invariant function
space $E=E(0,\infty )$ the fundamental function of $E$,
$\phi_{{E}}(\cdot )$, is given by
$$
\phi_{{E}}(t)=\Vert \chi_{{[0,t)}}\Vert_{{E}},\ t>0.
$$
%The dilation operator $D_s$ we define as in \cite{KPS} by $D_s
%f(t) = f(t/s)$, it is bounded on every r.i. space $E$ and
%$\norm{D_s}_{E\to E} \leq \max\set{1,s}$.
\subsection{Marcinkiewicz function and sequence spaces}\label{MFandSS}
Our main examples of fully symmetric function  and
sequence spaces are given in the following discussion.
Let $\Omega$ denote the set of concave functions $\psi : [0,\iny)
\to [0,\iny)$ such that $\lim_{t \to 0^+} \psi(t) = 0$ and
$\lim_{t \to \iny} \psi(t) = \iny$.
For $\psi \in \Omega$
define the weighted
mean function
$$
a\left( x,t\right) =\frac{1}{\psi \left(
t\right) } \int_{0}^{t}x^{\ast }\left( s\right) ds\quad t>0
$$
 and denote by
$M(\psi)$ the (Marcinkiewicz) space of measurable functions $x$ on
$[0,\iny)$ such that
\begin{equation}\label{WeieghedMeanF}
 \nm{x}_{M(\psi)} := \sup_{t
> 0} a\left( x,t\right) = \nm{a\left( x,\cdot\right)}_\iny
 < \iny.
\end{equation}
We assume in this paper that $\psi(t)=O(t)$ when $t\to 0$, which
is equivalent to the continuous embedding $M(\psi)\subseteq
L_\infty(J)$.
 The definition of the Marcinkiewicz sequence space
$m(\psi)$ of functions on $\mathbb N$
is similar,
$$
m(\psi)=\left\{x=\{x_n\}_{n=1}^\infty\ :\ \|x\|_{m(\psi)}:=\sup_{N\ge1}\frac{1}{\psi(N)}
\sum_{n=1}^Nx_n^*<\infty\right\}.
$$
\noindent {\bf Example (i)}. Introduce the following functions
$$
\psi_1(t)=\begin{cases}
t\cdot\log2,&0\le t\le 1\\
\log(1+t),& 1\le t <\infty
\end{cases},
$$
respectively (for $p>1$),
$$
\psi_p(t)=\begin{cases}
t,&0\le t \le 1\\
t^{1-\frac{1}{p}},& 1\le t <\infty
\end{cases}.
$$
The spaces
$
\clL^{(1,\infty)}$
%:=M(\log(1+t)) \cap L_\infty\quad
and $\clL^{(p,\infty)}$
%=M(t^{1-\frac 1p}) \cap L_\infty, \quad
%1<p<\infty
%$$
%Indeed, $\clL^{(1,\infty)}$
% (respectively, $\clL^{(p,\infty)}$, $p>1$) may be
are the Marcinkiewicz spaces $M(\psi_1)$ and $M(\psi_p)$
respectively. The norm given by  formula~\eqref{WeieghedMeanF} on
the space $\clL^{(p,\infty)}$ is denoted by $\|\cdot\|_{(p,\infty)}$,
$1\le p <\infty$.
%%%%%%%%%%%%%%%%%%%%%Kalton Sukochev II
%\subsection{Singular symmetric  functionals}
% on Marcinkiewicz spaces.}
%\label{SSFonMS}

\noindent{\bf Example (ii)}.
In \cite{N}, F. Nicola considers,
in connection with a class of
pseudo-differential operators, the Marcinkiwecz space $M(\psi)$,
with $\psi (t)=\log^2(t+1)$, $t>0$.

%\begin{defn}[cf. \cite{DPSS}, Definition 2.1]\label{SandRIf}
%A positive functional $f\in M(\psi)^*$ is said to be symmetric
%if $f(x) \leq f(y)$ for all $x,y \in
%M_+(\psi)$ such that $x\prec\prec y$.
%Such a functional is said to be supported at infinity (or
%singular) if $f(|x|) = 0$ for all $x \in M_1(\psi)$:= the closure
%of $L_1$ in $M(\psi)$ (equivalently, $f(x^*\chi_{[0,s]})=0$, for
%every $x\in M(\psi)$ and $s>0$ where $\chi_{[0,s]}$ is the
%indicator function of the interval $[0,s]$ ).
%\end{defn}
\subsection{Symmetric operator spaces and functionals.}
We now go from function spaces to the
setting of (noncommutative) spaces of operators.
Let $\clN $ be a semifinite von Neumann algebra on the separable
Hilbert space $\clH $, with a fixed faithful and normal
semifinite trace $\tau $.
%A closed densely defined
% operator $D$ affiliated to $\mathcal N$ is
%called $\tau$- measurable if for every $\epsilon >0$ there exists
%an orthogonal projection $p\in \clN$ such the $p({\clH })\subseteq$dom$(D)$ and $\tau(1-p)<\epsilon$. The set of all such is denoted $\clNt$.
We recall from ~\cite{FK,F} the notion of {\it generalized singular value
function}. Given a
self-adjoint operator $A$ in ${\mathcal N}$, we denote by
$E^A(\cdot )$ the spectral measure of $A$. %Now assume that $A$
%is$\tau$- measurable.
Then $E^{\vert A\vert}(B)\in \clN $ for all Borel sets
$B\subseteq{\mbR}$, and there exists $s>0$ such that $\tau
(E^{\vert A\vert}(s,\infty))<\infty$. For $t\geq 0$, we define
$$\mu_t(A)=\inf\{s\geq 0 : \tau (E^{\vert A\vert}(s,\infty))\leq t\}.$$
The function $\mu(A):[0,\infty)\to [0,\infty ]$ is called the {\it
generalized singular value function} (or decreasing rearrangement)
of $A$; note that $\mu_{(\cdot)}(A)\in L_\infty(J)$.

If we consider $\clN  =L_\infty([0,\infty),m)$, where $m$
denotes Lebesgue measure on $[0,\infty)$, as an abelian von
Neumann algebra acting via multiplication on the Hilbert space
${\clH }=L^2([0,\infty),m)$, with the trace given by
integration with respect to $m$, it is easy to see %that the set of
%all $\tau$-measurable operators affiliated with $\clN $
%consists of all measurable functions on $[0,\infty)$ which are
%bounded except on a set of finite measure, and
that the
generalized singular value function $\mu(f)$ is precisely the
decreasing rearrangement $f^*$.
If $\clN$ is all bounded operators (respectively, $\ell_\infty(\mbN)$)
and $\tau $ is the standard trace (respectively, the counting
measure on $\mbN$), then %it is not difficult to see that $\clNt
%=\clN$. In this case,
$A\in \clN$ is compact if and only if $\lim _{t\to \infty }\mu _t
(A)=0$; moreover,
$$
 \mu _n(A)=\mu _t(A), \quad t\in [n,n+1),\quad  n=0,1,2,\dots ,
$$
and the sequence $\{\mu _n(A)\}_{_{n=0}}^{\infty }$ is just the
sequence of eigenvalues of $\vert A\vert $ in non-increasing order
and counted according to multiplicity.

Given a semifinite von Neumann algebra $(\clN,\tau)$ and a fully
symmetric Banach function space $(E,\Vert\cdot\Vert_E)$ on
$([0,\infty),m)$, satisfying $E\subseteq L_\infty[0,\infty)$, we
define the corresponding non-commutative space $\ent $ by setting
$$E(\clN,\tau) = \{ A\in\clN : \mu (A)\in E\}.$$
The norm is $\Vert A\Vert _{_{E(\clN,\tau)}}:= \Vert \mu
(A)\Vert_E$, and the space $(E(\clN,\tau),\Vert\cdot\Vert
_{_{E(\clN,\tau)}})$ is called the (non-commutative) fully
symmetric operator space associated with $(\clN,\tau)$
corresponding to $(E,\Vert\cdot\Vert_E)$. We write $E(\clN,\tau)_+$
for the positive operators in $E(\clN,\tau)$. If
$\clN=\ell_\infty(\mbN)$, then the space $E(\clN,\tau)$ is simply
the (fully) symmetric sequence space $\ell_E$, which may be viewed
as the linear span in $E$ of the vectors $e_n=\chi_{_{[n-1,n)}}$,
$n\ge 1$ ({\it cf} \cite{LT}).

The spaces $M(\psi)\nt$ associated to Marcinkiewicz function spaces
are called operator Marcinkiewicz spaces  and
we mostly omit the symbol $\nt$
as this should not cause any confusion.
We use, for  the usual Schatten ideals in $\mathcal N$,
the notation  $L_p(\clN,\tau),\ p\geq 1$.

\begin{defn}
 A linear functional 
$\phi \in E(\clN,\tau )^*$ is called symmetric if $\phi$ 
is positive, (that is, $\phi (A)\geq 0$ whenever $0\leq A\in
E(\clN,\tau )$) and $\phi (A)\leq \phi (A')$ whenever 
$\mu(A)\prec\prec \mu (A')$. 
%whenever $A,A'\ge 0$ and $\mu (A)= \mu (A')$). 
A symmetric 
$\phi \in E(\clN,\tau )^*$ 
is called singular if it vanishes on all finite trace 
projections from $\clN$.
\end{defn}
The important examples of singular symmetric 
functionals that arise in noncommutative geometry
are the Dixmier traces which we describe in the next Section.
%\subsection{Singular symmetric functionals.} 
%Let
%$E=E(0,\infty)$ be a fully symmetric space. 
%By $E_+$ we denote the
%set of all  nonnegative functions from $E$ and by
%Let  $E^*_{sym}$ be the
%set of all symmetric functionals on $E$.
%A positive linear functional $\phi$ on $E$ is called \newn{normal} (or \newn{order continuous}),
% if from $f_n\downarrow 0$ it follows that $\phi(f_n)\downarrow 0.$
%A positive  $\phi\in E^*$ is said to be
%\newn{singular}, if from $0\leq \phi' \leq \phi$, $\phi'\in E^*,$
%$\phi'$ is order-continuous, it follows that $\phi' = 0.$
%\begin{prop} \cite{DPSS} (i)
%Every symmetric functional on $E$ can be uniquely decomposed into the sum of
%a normal functional and a singular symmetric functional.
%Moreover, the  normal functional is zero unless  $E\subseteq L^1[0,\infty).$\\
%(ii) Any singular symmetric functional can be uniquely decomposed into the
%sum of singular symmetric functionals, supported at zero and at infinity.\\
%(iii) The set of symmetric functionals forms a lattice.
%\end{prop}
For the discussion of these  we will need the following fact.
%Here $E\nt_{+}$
%will denote the set of all nonegative operators from $E\nt$.
%%%%%%%%%%%%%%%DPSS start
{
\def\p{\,\prec\kern-1.2ex \prec\,}
\def\ch{\raise 0.5ex \hbox{$\chi$}}
\def\m{\ ${\cal M}$}
\def\n{\ ${\cal N}$}
\def\nm{{\cal M}}
\def\mt{\ $\widetilde {\cal M}$}
\def\nmt{\widetilde {\cal M}}
\def\nn{{\cal N}}
\def\nnt{\widetilde {\cal N}}
\def\r{{\mbR}^+}
\def\h{H({\cal M})}
\def\g{G({\cal M})}
\def\up{M^{(p)}(E)}
\def\do{M_{(p)}(E)}
\def\doq{M_{(q)}(E)}
\def\em{E({\cal M},\tau )}
\def\mpc{$(L^p(\nm),{\cal M})$}
\def\mpcn{(L^p(\nm),{\cal M})}
\def\mp+{$L^p(\nm)+{\cal M}$}
\def\mp+n{L^p(\nm)+{\cal M}}
\def \cf {\rlap {$ \prec $}{$ \relbar$}}
\def \cfn {\rlap {\prec }{\relbar }}
\def \ycfx {$y {\rlap {\prec }{\relbar}} x$}
\def \ycfpx {$y {\rlap {\prec }{\relbar}} ^px$}

\begin{thm}[\cite{DPSS}]\label{gfasft7} {\it Let $\phi_0$
be a symmetric functional on $E$.
  If $\phi (A):=\phi _0(\mu (A))$, for all $A\geq 0$, $A\in E\nt$, then
  $\phi $ extends to a symmetric functional $0\leq \phi \in
  E(\clN,\tau )^*$.}
%Let $(\clN,\tau )$ be a
%semifinite
%AlanB
%von Neumann algebra without minimal projections, and let $E$ be a fully
%symmetric Banach function space on $[0,\infty )$. If $0 \leq \phi \in
%E(\clN,\tau )^*$ is a symmetric functional,
%AlanE
%then there exists a symmetric
%functional $0\leq \phi_0\in E^*$ such that $\phi (D)=\phi _0(\mu (D))$
%for all
%$0\leq D\in E(\clN,\tau )$.
\end{thm}

%Thus we may define singular functionals on  $E(\clN,\tau )$
%as those arising from singular functionals on the underlying function
%space $E$. 
%The correspondence between the sets
%$E^*_{sym}$ given in Theorems~\ref{gfasft7} and~\ref{gfasft8} also
%exists for the set of symmetric functionals $(\ell_E)^*_{sym}$ and
%$(E(\clH))^*_{sym}$. Furthermore, as the following result shows there exists
%a simple connection also between the sets $(\ell_E)^*_{sym}$ and $E^*_{sym}$.

%\begin{thm}[\cite{DPSS}] Let $E$ be a fully symmetric Banach function space on
%  $[0,\infty )$ and let $E(\clH)$ be the corresponding ideal of compact
%  operators on infinite-dimensional Hilbert space $\clH$.  If $0\leq \phi \in E\nt^*$ is a symmetric
%  functional, then there exists $\phi _0\in E^*_{sym}$ such that $\phi
%  (x)=\phi _0(\mu (x))$ for all $x\in E\nt_+$.
%\end{thm}

%}%

\section{Invariant states and Dixmier traces.} \label{PrelimDilInvState}
The construction of Dixmier traces $\tau_\omega$ depends crucially
on the choice of the ``invariant mean'' $\omega$.
%We
%denote by $\ell_\infty=\ell_\infty(\mbN)$ the Banach space of all
%bounded sequences of complex numbers.
%By a state on a unital
%$C^*$-algebra we mean a positive linear functional with value 1 on
%the unit of the algebra.
%We recall that a positive linear functional ${{\LL}}\in
%\ell_{\infty }(\mathbb{N})^{\ast }:=\ell_\infty^*$ is called a
%Banach limit if ${\LL}$ is translation invariant and
%${{\LL}}\left( \mathbf{1}\right) =1$ (here,
%$\mathbf{1}=(1,1,1\dots)$). A Banach limit ${\LL}$ satisfies in
%particular ${{\LL}}\left( \xi \right) =0$ for all $\xi \in c_{0}$
%($=$ all sequences from $\ell_\infty$ converging to zero).  Banach
%limits are in fact fundamental to a systematic study of Dixmier
%traces \cite{LSS}. 
Here we explain the invariance properties
we need for these invariant means
via the results summarized below (all of them are
proved using fixed
point theorems).
%We denote the collection of all Banach limits
%on $\ell_{\infty }$ by $BL(\mbN) $. Note that $\left\| \LL\right\|
%=1$ for all $ \LL\in BL(\mbN) $.

%\bigskip
%We recall that sequence $\xi =\left\{ \xi _{n}\right\}
%_{n=1}^{\infty}\in \ell_{\infty }$ is said to be {\it almost
%convergent}  to $\alpha\in \mbR$, denoted $\displaystyle
%F\text{-}\lim_{n\rightarrow \infty }\xi _{n}=\alpha $ if and only
%if ${{\LL}}\left( \xi \right) =\alpha$ for all ${{\LL}}\in BL(\mbN)$.
%The notion of an almost convergent sequence is due to G.G. Lorentz
%[5], who showed that the sequence $\left\{ \xi _{n}\right\}
%_{n=1}^{\infty}$ is almost convergent to $\alpha$ if and only if
%the equality $\displaystyle
%\lim_{p\rightarrow \infty }\frac{\xi _{n}+\xi _{n+1}+\cdots +\xi _{n+p-1}}{p}%
%=\alpha $ holds uniformly for $n=1,2,\dots\ $. We denote by $ac$
%(respectively, $ac_0$) the set of all almost convergent
%(respectively, all almost convergent to $0$) sequences from
%$\ell_\infty$. Clearly, $ac$ and $ac_0$ are closed subspaces in
%$\ell_\infty$.
We define the shift operator $\Shift\colon \ell_\infty\to \ell_\infty,$
the  Ces\`{a}ro operator $\Cesaro:\ell_\infty\to \ell_\infty$
and dilation operators $\Dilation_n:\ell_\infty \to \ell_\infty$ for
$n\in \mathbb{N}$ by the formulas
 %\label{ShiftHardyDilation}
$$  \Shift\Seque\x =(\x_2,\x_3,\x_4,\dots). $$
 $$ \Cesaro\Seque\x = (x_1,\frac {x_1+x_2}2,\frac {x_1+x_2+x_3}3,\dots), $$
 $$ \Dilation_n\Seque\x=(\underbrace{\x_1,\ldots,\x_1}_n,
\underbrace{\x_2,\ldots,\x_2}_n,\ldots),$$
for all $\x=\Seque \x \in \ell_\infty.$

\begin{thm}\cite{DPSSS2}
\label{InvState}
{\it There exists a state $\tilde\omega$ on $\ell_\infty$ such that for all $n\geq 1$
$$
  \tilde\omega\circ\Shift = \tilde\omega\circ\Cesaro = \tilde\omega\circ\Dilation_n =\tilde\omega.
$$}
\end{thm}

Now we consider analogous results for $L_\infty$. We let $\mbR^*_+$
denote the positive reals with multiplication as the group operation.
We define the isomorphism $L: L_\infty(\mbR)\to L_\infty(\mbR^*_+)$
by $L(f) = f\circ \log$. Next we define the Cesaro means (transforms)
on $L_\infty(\mbR)$ and $L_\infty(\mbR^*_+)$, respectively by:
$$H(f)(u)=\frac{1}{u}\int_0^u f(v) dv\quad \text{for}\quad f\in L_\infty(\mbR),\ u\in \mbR$$
and,
$$M(g)(t)=\frac{1}{\log t}\int_1^t g(s) \frac{ds}{s}\quad \text{for}\quad g\in
L_\infty(\mbR^*_+),\ t>0.$$
%We refer to $H$
%as the mean for the additive group $\mbR$.
A brief calculation yields for $g\in  L_\infty(\mbR^*_+)$,
$LHL^{-1}(g)(r)
%= \frac{1}{\log r}\int_0^{\log r} g(e^u)du
%= \frac{1}{\log r}\int_1^{r} g(v)\frac{dv}{v}
= M(g)(r),$
i.e  $L$ intertwines the two means.

%Analogues of the operators $\Shift, \Dilation_n$
%and $\Cesaro$ acting on $L_\infty(\mbR)$ and $L_\infty(\mbR_+^*)$
%are as follows.

\begin{defn}\label{familiesOfSelf-mapsDef}
Let $T_b$ denote translation by $b\in \mbR$, $D_a$ denote
dilation by $\frac{1}{a}\in \mbR^*_+$ and let $P^a$ denote exponentiation
by $a\in \mbR^*_+$. That is,
$$    T_b(f)(x) = f(x+b)\quad \text{for}\quad f\in L_\infty(\mbR),$$
$$    D_a(f)(x) = f\left(a^{-1}x\right)\quad \text{for}
\quad f\in L_\infty(\mbR),$$
 $$   P^a(f)(x) = f(x^a)\quad \text{for}\quad f\in L_\infty(\mbR^*_+).$$
\end{defn}

%The solvable transformation group generated by
%$\{D_a,P^a\ \vert\ \frac{1}{a}\in \mbR^*_+\}$ we call $G_2$.
%It plays a key role later on.
%Some of the basic relations between these $L_\infty$ spaces and their
%self-maps are provided for easy access by the following proposition,
%whose proof is similar to Lemma \ref{DilCesShiftLemma}.

%\begin{prop}[\cite{CPS2}] \label{DualityProp}
%$L_\infty(\mbR)$ together with the self-maps, $D_a$, $T_b$, and $H$
%($a>0,b\in\mbR$) is related to $L_\infty(\mbR^*_+)$ together with the
%self-maps, $P^a$, $D_a$, and $M$ ($a>0$) via the isomorphism
%$$L: L_\infty(\mbR)\to L_\infty(\mbR^*_+)$$ and the following identities:\\
%(1) $LD_{\frac{1}{a}}L^{-1} = P^a$ for $a>0$,\\
%(2) $LT_bL^{-1} = D_{(\exp(b))^{-1}}$ for $b\in\mbR$,
%(or $LT_{\log(a)}L^{-1} = D_a$ for
%$a>0$),
%\\
%(3) $LHL^{-1} = M$,\\
%(4) $D_aH = HD_a$ and $P^aM = MP^a$ for $a>0,$\\
%(5) $\liminfty{t} (HT_b-T_bH)f(t)=0$ for $f\in L_\infty(\mbR)$ and
 %   $b\in\mbR,$\\
%(6) $\liminfty{t} (MD_a-D_aM)f(t)=0$ for $f\in L_\infty(\mbR^*_+)$ and
%    $a>0.$
%\end{prop}

\begin{prop}[\cite{CPS2}]\label{cps2prop}
{\it If a continuous functional $\tilde\omega$ on $L_\infty(\mbR)$ is invariant under the Cesaro operator $\Cesaro,$
the shift operator $\Shift_a$ or the dilation operator $\Dilation_a$ then $\tilde\omega\circ L^{-1}$
is a continuous functional on  $L_\infty(\mbR^*_+)$ invariant under $\CesaroRPlus,$ the dilation operator $\Dilation_a$ or
$P^a$ respectively. Conversely, composition with $L$ converts an $\CesaroRPlus,$ $\Dilation_a$ or $P^a$ invariant
continuous functional on $L_\infty(\mbR^*_+)$
into an $\Cesaro,$ $\Shift_a$ or $\Dilation_a$ invariant continuous functional on  $L_\infty(\mbR)$.}
\end{prop}

We denote by $C_0(\mbR)$ (respectively, $C_0(\mbR^*_+)$) the continuous functions
on $\mbR$ (respectively, $\mbR^*_+$) vanishing at infinity (respectively at infinity and at zero).

\begin{thm}[\cite{CPS2}] \label{DualityThm}
{\it There exists a state $\tilde\omega$ on $L_\infty(\mbR)$
satisfying the following conditions:\\
(1) $\tilde\omega(C_0 (\mbR)) \equiv 0$.\\
(2) If $f$ is real-valued in $L_\infty(\mbR)$
then
$$ess \liminf\limits_{t\to\infty} f(t) \leq\tilde\omega(f)\leq ess \limsup\limits_{t\to\infty}
f(t).$$
(3) If the essential support of $f$ is compact then $\tilde\omega(f)=0.$\\
(4) For all $a>0$ and $c\in\mbR$
$
   \tilde\omega = \tilde\omega\circ T_c = \tilde\omega\circ D_a = \tilde\omega\circ H.
$}
\end{thm}
Combining Theorem~\ref{DualityThm} and Proposition~\ref{cps2prop}, we
obtain
\begin{cor} \label{DualityCor}
{\it There exists a state $\omega$ on $L_\infty(\mbR^*_+)$
satisfying the following conditions:\\
(1) $\omega(C_0(\mbR^*_+)) \equiv 0$.\\
(2) If $f$ is real-valued in $L_\infty(\mbR^*_+)$
then
$$ess\liminf\limits_{t\to\infty} f(t) \leq\omega(f)\leq ess \limsup\limits_{t\to\infty} f(t).$$
(3) If the essential support of $f$ is compact then $\omega(f)=0.$\\
(4) For all $a,c>0$
$
   \omega = \omega\circ D_c = \omega\circ P^a = \omega\circ M.
$}
\end{cor}
\begin{rems}\label{pairs_of_functionals} In the sequel we will consider pairs of functionals
$\tilde\omega$ on $L_\infty(\mbR)$, $\omega\in L_\infty(\mbR^*_+)$
related by  $\tilde\omega\circ L^{-1}=\omega$.
\end{rems}

If $\omega $ is a state on $\ell_\infty$ (respectively, on
$L_\infty(\mbR)$, $L_\infty(\mbR_+^*)$), then we denote its value
on the element $\{x_i\}_{i=1}^\infty$ (respectively, $f\in
L_\infty(\mbR)$, $L_\infty(\mbR_+^*)$) by $\omega
-\lim_{i\to\infty} x_i$ (respectively, $\omega
-\lim_{t\to\infty}f(t)$). We saw in  Theorems~\ref{InvState},
\ref{DualityThm} and Corollary~\ref{DualityCor} states on
$\ell_\infty$, $L_\infty(\mbR)$, and $L_\infty(\mbR_+^*)$
invariant under various (group) actions. Alain Connes
in~\cite{Co4} suggested working with the set of states
on~$L_\infty(\mbR_+^*)$, which is larger then the set
$$\{\omega:\ \ \text{$\omega$ is an $M$-invariant state
    on~$L_\infty(\mbR_+^*)$}\}$$ namely
$$ CD(\mbR_+^*) :=
\{\tilde \omega = \gamma
\circ M\ :\ \text{$\gamma$ is an arbitrary singular state
  on~$C_b[0, \infty)$}\}. $$
These states are automatically dilation invariant.
In this paper, we find that for the zeta function
asymptotics it suffices to consider  states that are
$D_2$ and $P^\alpha$
invariant for all $\alpha>1$.

In Section 5 we need
 a smaller set of states, namely a subset of
$$
\{\omega \in L_\infty(\mbR_+^*)^*:\ \ \text{$\omega$ is an
  $M$-invariant and $P^a$-invariant state on~$L_\infty(\mbR_+^*)$, $a
  >0$}\}.
$$
This subset consists of  states whose existence is guaranteed by
Corollary \ref{DualityCor}. We refer to any state satisfying the
conditions (1) to (4) of Corollary \ref{DualityCor} as a {\it
DPM state} (in \cite{CPS2} we used the vaguer term `maximally invariant').
We now recall the construction of Dixmier traces for the compact
operators.
\begin{defn}\label{dixtrd}
Let $\omega$ be a $D_2$-invariant state on $\ell_\infty$.
The associated Dixmier trace of $T \in \clL_+^{(1,\infty)}(\clH)$ is the number
$$
\tau_\omega(T) := \limo{N}\ILogN \sums{n=1}^N \mu_n(T).$$
\end{defn}
Notice that in this definition we have chosen $\omega$ to satisfy only the
  dilation invariance assumption
even though Dixmier \cite{Dix66CR} originally imposed on $\omega$ the
  assumption of dilation and translation invariance.

Definition~\ref{dixtrd} extends to the Marcinkiewicz spaces
$M(\psi)\nt$. Fix an arbitrary
$D_2$-invariant state $\omega$ on $L_\infty(\mbR_+^*)$. Then the
state $\omega$ is $D_{2^n}$-invariant, $n\in{\mathbb Z}$} and a
simple argument shows that it also satisfies conditions (1)--(3)
of Corollary~\ref{DualityCor}. For the remainder of the paper, let
$\psi\in \Omega$ satisfy
\begin{equation} \label{psicon}
\lim_{t\to\infty}\frac{\psi(2t)}{\psi(t)}=1.
\end{equation}
This condition is  sufficient for the existence of singular traces
or singular  symmetric functionals on the corresponding fully
symmetric operator spaces \cite{DPSS}. Indeed, setting
\begin{equation}
\label{Troeq} \tau_\omega(x):=\omega\text{-}
\lim_{t\to\infty}a(x,t),\quad 0\le x\in M(\psi)\nt
\end{equation}
(see the details in~\cite[p.~51]{DPSS}), we obtain an additive
homogeneous functional on $M(\psi)\nt_+$, which extends to a
symmetric functional on $M(\psi)\nt$ by linearity. The proof of
linearity of $\tau_\omega$ in~\cite[p.~51]{DPSS} is based on the
assumption that $\omega$ is $D_{\frac{1}{2}}$-invariant which is
equivalent to $D_2$-invariance (see above).
\section{The Dixmier trace on Marcinkiewicz  operator spaces }
\subsection {Preliminaries}
In this subsection we generalize and strengthen some results from
\cite{CPS2}.

\begin{lemma} \label{basic}
{\it For every $\psi\in \Omega$ satisfying \eqref{psicon} and
every $1>\alpha>0$, there is $C=C(\alpha)$ such that
$\psi(t)<Ct^\alpha,\ t>0$.}
\end{lemma}

\begin{proof} Let $0<\alpha$ and let $Q>0$ be so large that  for $t>Q$
$$
\frac{\psi(2t)}{\psi(t)}< 2^{\alpha}.
$$
 There is $C>1$ so large that $\psi(t)\le Ct^{\alpha}$ for all $t<Q$.
Suppose there is a first $Q_0\ge Q$ for which $\psi(Q_0)=CQ_0^{\alpha}$. Then
$$
\frac{\psi(Q_0)}{\psi(Q_0/2)}\ge\frac{CQ_0^{\alpha}}{C(Q_0/2)^{\alpha}}
=2^\alpha,
$$
which is a contradiction.
Consequently, $\psi(t)<Ct^{\alpha}$ for all $t>0$. %Since $t^{\alpha'}=o(t^\alpha$ we have $\psi(t)<t^{\alpa}$ for all sufficiently large $t$.
\end{proof}

\noindent Recall that for any $\tau$-measurable operator $T$, the distribution
function of $T$ is defined by
$$
\lambda _t(T):=\tau (\chi_{(t,\infty)}(|T|)),\quad t>0,
$$
where $\chi_{(t,\infty)}(|T|)$ is the spectral projection of $|T|$
corresponding to the interval $(t,\infty)$ (see \cite{FK}). By
Proposition 2.2 of \cite{FK},
$$\mu_s(T)=\inf\{t\ge 0\ :\ \lambda_t(T)\leq s\}.
$$
We infer that for any  $\tau$-measurable operator $T$, the distribution
function $\lambda _{(\cdot)}(T)$ coincides with the (classical) distribution
function of $\mu_{(\cdot)}(T)$. From this formula and the fact that $\lambda$
is right-continuous, we can easily see that for $t>0$, $s>0$
$$s\geq\lambda_t\Longleftrightarrow \mu_s\leq t.$$
Or equivalently,
$$s < \lambda_t \Longleftrightarrow \mu_s > t.$$
Using Remark 3.3 of \cite{FK} this implies that:
\begin{equation}\label{six}
\int_0^{\lambda _t} \mu_s(T)ds=\int_{[0,\lambda _t)} \mu_s(T)ds=
\tau (|T|\chi_{(t,\infty)}(|T|)),\quad t>0.
\end{equation}

%Let $\lambda_{(.)}(T)$ be the spectral
% distribution function of a positive operator $T\in M(\psi)$.

\begin{lemma}
{\it For $T\in M(\psi)$ $T\geq 0$
and any $\beta>1$ there is a $C=C(\beta)$ such that
$
\lambda_{1/t}(T)<Ct^\beta
$
for every $t>0$. }
\end{lemma}
\begin{proof} Let $\alpha=1-1/\beta$ and $\lambda_{1/t}(T)=a$ .
Hence $\mu_{(a-0)}(T)\ge 1/t$.
Then by Lemma \ref{basic} there is $C_1>0$ such that
$$
\|T\|_\psi=\sup_{0<h<\infty}\frac{\int_0^h \mu_s(T)ds}{\psi(t)}\ge
\frac{\int_0^a \mu_{(a-0)}(T)ds}{\psi(a)}= \frac
{a\mu_{(a-0)}(T)}{\psi(a)}\ge
\frac{a(1/t)}{C_1a^\alpha}=a^{1-\alpha}/(C_1t).
$$
Consequently
$$
 \lambda_{1/t}(T)=a<(C_1\|T\|_\psi t)^{1/(1-\alpha)}=Ct^\beta.
$$
\end{proof}

\noindent{\bf Remark}.
 Since $\beta>1$ could be arbitrary, it is obvious that
the constant $C$ could be  replaced by 1 if $t$ is sufficiently
large.

In the sequel we will suppose that $\psi$ possesses the following property
\begin{equation}\label{A}
A(\beta)=sup_{t>0}\frac{\psi(t^\beta)}{\psi(t)}\to 1,
\mbox{ if } \beta\downarrow 1.
\end{equation}
Observe that if $\psi(t)=\log(1+t)^\gamma,\ \gamma>0$, then condition
(\ref{A}) is satisfied.

\begin{prop} (cf. \cite[Proposition 2.4]{CPS2})\label{prop2.4}
{\it For $T\in{\mathcal M}(\psi)$ positive let $\omega$ be  $D_2$
and $P^\alpha$-invariant, $\alpha>1$ state on $L^\infty(
\mathbb{R}^*_+)$. Then
$$
\tau_\omega(T)=
\omega-\lim_{t\to\infty}\frac{1}{\psi(t))}\int_0^t \mu_s(T) ds
=\omega-\lim_{t\to\infty}\frac{1}{\psi(t)}
\tau(T\chi_{(\frac{1}{t},\infty)}(T))
$$
%$$
%=\omega-\lim_{t\to\infty}\frac{1}{\psi(t)}\int_0^{t^\beta}\mu_s(T) ds
%$$
and if one of the $\omega-$limits is a true limit then so is the
other. }
\end{prop}

\begin{proof} We first note that
$$
\int_0^t\mu_s(T)ds\leq \int_0^{\lambda _{\frac{1}{t}}(T)} \mu_s(T)ds+1,
\quad t>0.
$$
Indeed, the inequality above holds trivially if
$t\leq\lambda _{\frac{1}{t}}(T)$. If $t> \lambda _{\frac{1}{t}}(T)$, then
$$
\int_0^t\mu_s(T)ds=\int_0^{\lambda _{\frac{1}{t}}(T)} \mu_s(T)ds+\int_
{\lambda _{\frac{1}{t}}(T)}^t \mu_s(T)ds.
$$
Now $s> \lambda _{\frac{1}{t}}(T)$ implies that $\mu_s(T)\leq \frac{1}{t}$
so we have
$$
\int_0^t\mu_s(T)ds\leq \int_0^{\lambda _{\frac{1}{t}}(T)} \mu_s(T)ds+\frac{1}{t}(t-\lambda _{\frac{1}{t}}(T))\leq \int_0^{\lambda _{\frac{1}{t}}(T)} \mu_s(T)ds+1.
$$
Using this observation and lemma and remark  above we see that for
%$C>\Vert T\Vert _{{\mathcal L}^{(1,\infty)}}$ and any fixed
$\alpha >1$  eventually
$$
\int_0^t\mu_s(T)ds \leq \int_0^{\lambda _{\frac{1}{t}}(T)} \mu_s(T)ds+1 %\int_0^{Ct\log t} \mu_s(T) ds +1
\leq\int_0^{t^\alpha}
\mu_s(T) ds+1
$$
and so eventually
$$
\frac{1}{\psi(t)}\int_0^t\mu_s(T)ds
\leq\frac{1}{\psi(t)}( \int_0^{\lambda _{\frac{1}{t}}(T)} \mu_s(T)ds+1)
\leq\frac{1}{\psi(t)}(\int_0^{t^\alpha} \mu_s(T) ds +1)
$$
$$
\leq\frac{\psi(t^\alpha)}{\psi(t)\psi(t^\alpha)}
(\int_0^{t^\alpha} \mu_s(T) ds+1).
$$
Taking the $\omega$-limit we get
$$
\tau_\omega(T) \leq\omega-\lim_{t\to\infty}\frac{1}{\psi(t)}
\int_0^{\lambda _{\frac{1}{t}}(T)} \mu_s(T)ds
\leq\omega-\lim_{t\to\infty}\frac{1}{\psi(t)}\int_0^{t^\alpha} \mu_s(T) ds
$$
$$
\leq\omega-\lim_{t\to\infty}\frac{A(\alpha)}{\psi(t^\alpha)}
\int_0^{t^\alpha} \mu_s(T) ds=A(\alpha)\tau_\omega(T)
$$
where the last line uses $P^\alpha,\ \alpha>1,$ invariance. %Remark 3.3 of \cite{FK}
%implies that:
%$$ \int_0^{\lambda _u} \mu_s(T)ds=\int_{[0,\lambda _u)} \mu_s(T)ds=
%\tau (|T|\chi_{(u,\infty)}(|T|)),\quad u>0.$$
Due to equality \eqref{six} and  since the previous inequalities
hold for all $\alpha>1$ and by assumption \eqref {A} we have
$A(\alpha)\to 1$ we get the conclusion of the proposition for
$\omega$-limits.

To see the last assertion of the Proposition suppose that
$\lim_{t\to\infty}\frac{1}{\psi(t)}\int_0^t\mu_s(T)ds=B$
then by the above argument for any $\epsilon>0$ and sufficiently large $t>0$ we get
$$
B-\epsilon\leq\frac{1}{\psi(t)}
\tau(T\chi_{(\frac{1}{t},\infty)}(T))
%\leq\limsup_{t\to\infty}\frac{1}{\log(1+t)}
%\tau(T\chi_{(\frac{1}{t},\infty)}(T))
\leq A(\alpha) (B+\epsilon)
$$
for all $\alpha>1$ and since $A(\alpha)\to 1$,
$\lim_{t\to\infty}\frac{1}{\psi(t)}
\tau(T\chi_{(\frac{1}{t},\infty)}(T))=B$.
On the other hand if the limit
$\lim_{t\to\infty}\frac{1}{\psi(t)}\tau(T\chi_{(\frac{1}{t},\infty)}(T))$
exists and equals $B$ say then
$$
\limsup_{t\to\infty}\frac{1}{\psi(t)}\int_0^t \mu_s(T) ds\leq B\leq
A(\alpha)\liminf_{t\to\infty}\frac{1}{\psi(t)}\int_0^t \mu_s(T) ds
$$
for all $\alpha>1$ and so
$
\lim_{t\to\infty}\frac{1}{\psi(t)}\int_0^t \mu_s(T) ds=B
$
as well. %The remaining claims follow similarly.
\end{proof}

\begin{cor}\label{after_prop2.4} {\it Under the conditions of
the preceding Proposition the expression
$$
\omega-\lim_{t\to\infty}\frac{1}{\psi(t)}
\tau(T\chi_{(\frac{1}{t},\infty)}(T))
$$
can be replaced by
$$
\omega-\lim_{t\to\infty}\frac{1}{\psi(t)}
\tau(T\chi_{(\frac{1}{t},1)}(T)).
$$
If the real limit exists then the prefix $\omega$ may be removed.}
\end{cor}
The proof is immediate since $\psi(\infty)=\infty$ and the
difference of these limits is
$$
\lim_{t\to\infty}\frac{1}{\psi(t)}\tau(T\chi_{(1,\infty)}(T))=\lim_{t\to\infty}\frac{1}{\psi(t)}\int_0^{\lambda_1(T)}\mu_s(T)ds=0.
$$
\subsection{An  alternative description of $\clL ^{(1,\infty)}$.}
% There is an elementary estimate on the singular values of
%operators $T$ in $\clL ^{(1,\infty)}$, which follows immediately
%from the submajoriation $\mu_s(T)\prec\prec K/(1+s)$, $K>0$ and
%well-known result of Hardy-Littlewood-P\' olya (see e.g. \cite{F},
%Lemma 4.1).
%
%\begin{lemma}\label{2.1}
%{\it For $T\in\clL ^{(1,\infty)}$ positive there is a constant
%$K>0$ such that for each $p\geq 1$, }
%$$\int_0^t\mu_s(T)^pds \leq K^p\int_0^t\frac{1}{(s+1)^p}ds.$$
%\end{lemma}
%\begin{proof}  By \cite{FK},
%Lemma 2.5 (iv),  for all $0\leq T\in {\mathcal N}$ and
%all continuous increasing functions $f$ on $[0,\infty)$ with $f(0)\ge 0$,
%we have $\mu_s(f(T))=f(\mu_s(T))$ for all $s>0$.
%Combining this fact with well-known result of Hardy-Littlewood-P\' olya
%(see e.g. \cite{F}, Lemma 4.1), we see that $T_1\prec\prec T_2$,
%$0\leq T_1,T_2\in {\mathcal N}$ implies $T_1^p\prec\prec T_2^p$ for
%all $p\in (1,\infty)$. Now, by definition of ${\mathcal L}^{(1,\infty)}$
%the singular values of $T$ satisfy $\int_0^t\mu_s(T)ds=O(\log t)$ so that
%for some $K>0$,
%$$\int_0^t\mu_s(T) ds \leq K \int_0^t\frac{1}{(s+1)}ds,\quad \forall t>0.$$
%In other words $\mu_s(T)\prec\prec K/(1+s)$ and the assertion of lemma
%follows immediately.
%\end{proof}
The zeta function of a positive compact operator $T$ is given by
$\zeta(s)=\tau(T^s)$ for real positive $s$ on the assumption that
there exists some $s_0$ for which the trace is finite. Note that
it is then true that $\tau(T^s)<\infty$ for all $s>s_0$.
 In this subsection we will always
assume $\tau(T^s)<\infty $ for all $s>1$ and we  are interested in
the asymptotic behavior of $\zeta(s)$ as $s\to 1$.
% If $T\in \clL
%^{(1,\infty)}$, by Lemma \ref{2.1} we have for some $K>0$ and all
%$s>1$
% $\tau(T^s) \leq \frac{K^s}{s-1}.$
%Thus for $\tilde \omega\in L_\infty(\mbR)^*$
%\begin{equation}\label{(3.1)}
%\tilde\omega-\lim_{r\to \infty}\frac{1}{r}\tau(T^{1+\frac{1}{r}})
%\quad \end{equation}
% exists.

%At first look at Lemma 3.3(i) [CPS2] (think, $b=A$ there). The
%assumption (the one which Alan wants to have from the outset)
%$$\lim_{s\to 0} s\tau(A(1+D^2)^{-(1+s)/2})=C<\infty $$
% implies that
%$$\sup_{s>0}\tau((A^{1/2}(1+D^2)^{-1/2} A^{1/2})^{1+s})<\infty$$
%(via [CPS2] Lemma 3.3) and that the limit $\lim_{s\to 0} s
%\tau((A^{1/2}(1+D^2)^{-1/2} A^{1/2})^{1+s})$  exists. By
%Proposition \ref{newresult}, the existence of the latter limit
%yields $$A^{1/2}(1+D^2)^{-1/2} A^{1/2}\in \clL ^{(1,\infty)}.$$

%Take $a \in B(H)^+$, $\dim Ker (a )= \dim H$, $a \in \clL
%^{(p,\infty)}$, $\sqrt{a} \notin \clL ^{(p,\infty)}$.

%We now define an orthogonal projection $y$ in $H\oplus H$ by
%setting $$y=[a_{ij}]_{i,j=1}^2,\quad a_{11}=a,\quad
%a_{12}=\sqrt{a-a^2},\quad a_{22}=1_H -a.$$ We also define an
%orthogonal projection $x$ in $H\oplus H$ by setting,
%  $$x=[b_{ij}]_{i,j=1}^2,\quad
%  b_{11}=1,\quad b_{12}=b_{22}=0.$$
%We clearly have $a=xyx\in \clL ^{(p,\infty)}$. On the other hand
%for $s$-numbers of the operator $xy$, we have
%$s_n(xy)=s_n(\sqrt{a})$ $\forall n \in {\mathbb N}$, due to the fact
%$a=xyx=(xy)(xy)^*$. Thus, $xy\notin \clL ^{(p,\infty)}$.

Let us define the space
$$
{\mathcal Z}_1=\{T\in {\mathcal N}:\ \|T\|_{{\mathcal
Z}_1}=\limsup_{p\downarrow 1}(p-1)\tau(|T|^p)<\infty\}.
$$
Since we also have the other equivalent definition
$$\|T\|_{{\mathcal Z}_1}=\limsup_{p\downarrow 1}(p-1)(\int_0^\infty \mu_t(|T|)^pdt)^{1/p}=\limsup_{p\downarrow1}(p-1)\|T\|_{L_p}$$
(recall that we use the notation $L_p$ for the Schatten ideals in
$(\clN,\tau)$) the ordinary properties of the semi-norm for
$\|\cdot\|_{{\mathcal Z}_1}$ are immediate.
\begin{thm}\label{newresult}
(i) {\it Let $T\geq 0$, $T\in \mathcal N$
and
$\limsup_{s\to 0}s\tau(T^{1+s})=C<\infty,$ then}
%$x\in{\cal L}^{1,\infty}$ and
$$
\limsup_{u\to\infty}\frac 1{\ln u}\int_0^u \mu_{t}(T)dt\le Ce.
$$
(ii) {\it The spaces ${\mathcal Z}_1$ and ${\mathcal
L}^{1,\infty}$ coincide. Moreover, if  $\mathcal {N}$ is a type
$I$ factor with the standard trace, or else $\mathcal {N}$ is
semifinite and the trace is non-atomic then denoting by ${\mathcal
L}^{1,\infty}_0$ the closure of $L_1(\mathcal{N},\tau)$ in
${\mathcal L}^{1,\infty}$, we have for any $T\in {\mathcal C}_1$ }
$$
\mbox{dist}_{{\mathcal L}^{1,\infty}}(T,{\mathcal L}^{1,\infty}_0)
=\limsup_{u\to\infty}\frac 1{\ln u}\int_0^u \mu_t(T)dt\le
e\|T\|_{\mathcal{Z}_1}
$$
{\it and} $\|T\|_{{\mathcal Z}_1}\le \|T\|_{1,\infty}$.
\end{thm}

\begin{proof}
\noindent
(i) By assumption for every $\epsilon>0$ there is
an $s_0>0$ such that for all $s\in [0,s_0]$
\begin{equation}\label{(7)}
\quad s\int_0^\infty \mu_t(T)^{1+s}dt\le C+\epsilon.
\end{equation}
Then, for $u\ge 1$ according to H\"older's inequality and (\ref{(7)})
we have
$$
\int_0^u\mu_t(T)dt\le
\left(\int_0^u\mu_t(T)^{1+s}dt\right)^{\frac1{1+s}}
\left(\int_0^u 1^{\frac{1+s}s}dt\right)^{\frac s{1+s}}\le
$$$$
\left(\frac s{s}\int_0^\infty
\mu_t(T)^{1+s}dt\right)^{\frac1{1+s}}u^{\frac s{1+s}}\le
((C+\epsilon)/s)^{\frac1{1+s}}u^{\frac s{1+s}}\le (C+\epsilon)\frac 1su^s.
$$
Set $u_0={e^{1/{s}_0}}$ and for $u>u_0$ set $s=1/\ln u(<s_0)$. Then $u=e^{\ln u}$ and by the previous inequality
$$
\int_0^{u}\mu_t(T)dt\le (C+\epsilon)\frac 1su^s= (C+\epsilon)\frac
{e^{\ln u\frac 1{\ln u}}}{\frac 1{\ln u}}=(C+\epsilon)e\ln u.
$$
That is we have the inequality
$$
\frac1{\ln u}\int_0^u \mu_t(T)dt\le (C+\epsilon)e \mbox{ for } u>u_0.
$$
Since $$\|T\|_{{\mathcal L}^{1,\infty}}
=\sup_{1\le u\le\infty}\frac 1{\ln(1+u)}\int_0^u\mu_t(T)dt$$
we conclude that $T\in{\mathcal L}^{1,\infty}$.
Moreover, since $\epsilon>0$ is arbitrary
$$
\limsup_{u\to\infty}\frac 1{\ln u}\int_0^u \mu_t(T)dt\le eC.
$$
Hence (i) and the embedding $ {\mathcal Z}_1\subset{\mathcal
L}^{1,\infty}$ are established.

The equality $ \mbox{dist}_{{\mathcal L}^{1,\infty}}(T,{\mathcal
L}^{1,\infty}_0) =\limsup_{u\to\infty}\frac 1{\ln u}\int_0^u
\mu_t(T)dt$ is well-known in the special case when the algebra
$\mathcal {N}$ is commutative (see e.g. \cite[Proposition
2.1]{DPSSS1} and references therein). The general case follows
from this special case, due to the combination of the following
facts. Firstly, the inequality $\mu(x)-\mu(y)\prec\prec \mu(x-y)$
(see \cite{DDP}) together with the fact that ${\mathcal
L}^{1,\infty}$ is fully symmetric yields the inequality $
\mbox{dist}_{{\mathcal L}^{1,\infty}}(T,{\mathcal
L}^{1,\infty}_0)\ge  \mbox{dist}_{{\mathcal
L}^{1,\infty}}(\mu(T),{\mathcal L}^{1,\infty}_0(0,\infty))$ or $
\mbox{dist}_{{\mathcal L}^{1,\infty}}(T,{\mathcal
L}^{1,\infty}_0)\ge  \mbox{dist}_{{\mathcal
L}^{1,\infty}}(\mu(T),{\mathcal L}^{1,\infty}_0(\mathbf{N}))$,
depending whether $\mathcal{N}$ is of type $II$ or $I$. Secondly,
fix an arbitrary $T\in {\mathcal L}^{1,\infty}(\mathcal{N})$. Due
to \cite{CS}, there exists a rearrangement-preserving (and thus,
isometric) embedding $\phi_T$ of ${\mathcal
L}^{1,\infty}(0,\infty)$ (respectively, ${\mathcal
L}^{1,\infty}(\mathbf{N})$ in the type $I$ setting) into
${\mathcal L}^{1,\infty}(\mathcal{N})$ such that
$\phi_T(\mu(T))=T$. This observation shows that $
\mbox{dist}_{{\mathcal L}^{1,\infty}}(T,{\mathcal
L}^{1,\infty}_0)\le  \mbox{dist}_{{\mathcal
L}^{1,\infty}}(\mu(T),{\mathcal L}^{1,\infty}_0(0,\infty))$.

The argument above also proves the equality and the first
inequality in (ii).

To complete the proof of (ii), let us take an arbitrary $T\in
{\mathcal L}^{1,\infty}$ and note that by the definition of the
norm in the Marcinkiewicz space ${\mathcal L}^{1,\infty}$ we have
$x\prec\prec \|T\|_{1,\infty}/(1+t)$. Since the spaces
$L_p(\mathcal{N},\tau)$, $ 1\le p\le\infty$, are fully symmetric
operator spaces we have
$$
\|T\|_p\le \|T\|_{1,\infty}\|1/(1+t)\|_p,\ p>1.
$$
Taking the $p$-th power we get
$$
\int_0^\infty \mu_t(T)^p\,dt\le
\|T\|_{1,\infty}^p\int_0^\infty1/(1+t)^pdt=\|T\|_{1,\infty}^p\frac1{p-1}.
$$
If now $p\downarrow 1$ we conclude that
$$
\|T\|_{{\mathcal Z}_1}=\limsup_{p\downarrow 1}(p-1)\int_0^\infty
\mu_t(T)^p \,dt\le \|T\|_{1,\infty}.
$$
Hence, ${\mathcal L}^{1,\infty}\subset {\mathcal Z}_1$. Due to the
first part of the proof we infer that the spaces ${\mathcal Z}_1$
and ${\mathcal L}^{1,\infty}$ are coincident.
\end{proof}

\begin{cor}\label{bruno}
{\it Let $T\in\mathcal N$ be positive with
$\tau(T^s)<\infty $ for all $s>1$.
 If $\lim_{r\to\infty}\frac{1}{r}\tau(T^{1+\frac{1}{r}})$
exists then $T\in\clL ^{(1,\infty)}$. }
\end{cor}
\subsection{The case $p>1$.}
 Our approach above to the
study of ${\mathcal Z}_1$ allows us to generalize immediately. Let
us define a class of spaces ${\mathcal Z}_q,\ q\ge 1$ by:
$$
{\mathcal Z}_q=\{T\in {\mathcal N}_+:\ \|T\|_{{\mathcal Z}_q}=
\limsup_{p\downarrow q}((p-q)\tau(T^p))^{1/p}<\infty\}.
$$
Setting $r=1+\frac{p-q}{q}=\frac pq$, we have
$$
\| T\|_{{\mathcal Z}_q}= \limsup_{p\downarrow
q}((p-q)\tau(T^{q(1+(p-q)/q)}))^{1/p}= (q\limsup_{p\downarrow
q}(p-q)/q\tau((T^q)^{(1+(p-q)/q)}))^{1/p}
$$$$
=q^{1/q}(\limsup_{r\downarrow 1}((r-1)\tau
((T^q)^r))^{1/(qr)}=(q\|T^q\|_{{\mathcal Z}_1})^{1/q}.
$$
Now it is clear that $T\in {\mathcal Z}_q$ if and only if $T^q\in
{\mathcal Z}_1$ and $\|T\|_{{\mathcal Z}_q}=(q\|T^q\|_{{\mathcal
Z}_1})^{1/q}$.

We now state a few consequences of Theorem \ref{newresult}.
%We make one remark which may help the reader understand what we are doing.
%The classical  $p$-convexification procedure for an arbitrary
%Banach lattice $X$ appears in \cite[Section 1.d]{LT2}.
%It is simply a direct generalization of the procedure of
%defining $L_p$-spaces from an $L_1$-space, and is applicable to any
%Banach lattice and/or r.i. operator space. Exactly the same
%process which is used to define $L_p$-spaces from
%a given $L_1$ (say on a measure space in the classical measure
%theory; or on semifinite von Neumann algebra in non-commutative
%integration theory) leads to ${\mathcal Z}_p$ if you start with
%${\mathcal L}^{1,\infty}$. In this point of view ${\mathcal
%L}^{1,\infty}$ is an "$L_1$-space", only defined through a Dixmier
%trace (that is $x\in {\mathcal L}^{1,\infty}$ if and only if
%$\tau_\omega(|x|)<\infty$) and not the classical trace; and so one
%can regard ${\mathcal Z}_p$ as analogues of an $L_p$-space
%with respect to a given Dixmier trace.
The classical  $p$-convexification procedure for an arbitrary
Banach lattice $X$ is described in \cite[Section 1.d]{LT2} and is
sometimes termed power norm transformation. It is simply a direct
generalization of the procedure of defining $L_p$-spaces from an
$L_1$-space.

The
proof of the first corollary below is immediate.
\begin{cor}\label{OrM}
(i) {\it There is a more convenient equivalent formula for the
semi-norm $\|\cdot\|_{{\mathcal Z}_q}$ namely}
$$
\|T\|_{{\mathcal Z}_q}^+=\|T^q\|_{{\mathcal Z}_1}^{1/q},\ q\ge 1.
$$
(ii) {\it The  space ${\mathcal Z}_q$ coincides  as a set with the
 $q$-convexification  of the operator space ${\mathcal
L}^{1,\infty}$ :
$$
{\mathcal L}_q^{1,\infty}=\{T\in {\mathcal N}_+:\
\|T\|^q_{1,\infty}=\sup_{1<u<\infty}\left(\frac{\int_0^u
\mu_t(T)^qdt}{\log(1+u)}\right)^{1/q}<\infty\}.
$$
If
 $\mathcal
{N}$ is a type $I$ factor with the standard trace, or else
$\mathcal {N}$ is semifinite and the trace is non-atomic then the
semi-norms $\|\cdot\|_{{\mathcal Z}_q}$ and} dist$_{{\mathcal
L}_q^{1,\infty}}(\cdot,{\mathcal L}^{1,\infty}_{q,0})$ {\it are
equivalent. Here,  ${\mathcal L}^{1,\infty}_{q,0}$ is the closure
of $L_1(\mathcal{N},\tau)$ in ${\mathcal L}_q^{1,\infty}$. }
\end{cor}

\begin{cor} {\it (i) An element $T\in {\mathcal Z}_p,\ p\ge 1,$
iff $T^p\in {\mathcal L}^{1,\infty}$.
Moreover
\begin{equation}\label{(+)}
\frac 1r\int_0^\infty
\mu_t(T)^{p+1/r}dt=\frac1r\tau(T^{p+1/r})=p\frac1{pr}\tau({T^p}^{(1+1/pr)}).
\end{equation}
and for $r>0$
the expression in (\ref{(+)}) belongs to $L^\infty(\mathbb{R}^*_+)$.\\
(ii) If $T\in {\mathcal L}^{p,\infty}$ then $T\in {\mathcal Z}_p$.\\
(iii) If  $T$ is a positive in $\mathcal N$ such
that $\lim_{r\to\infty}\frac{1}{r}\tau(T^{p+\frac{1}{r}})$ exists,
then $T\in{\mathcal Z}_p$. }
\end{cor}

\begin{proof}
The first statement is immediate from earlier results.
To prove (ii) we remind
the reader that
$T\in {\mathcal L}^{p,\infty}$
iff $\mu_t(T)\le C\min(1,t^{-1/p})$ for some $C<\infty$. Then as $r\to\infty$
$$
\frac 1r\int_0^\infty {\mu_t(T)}^{p+1/r}dt\le C\frac 1r(1+\int_1^\infty t^{-1-1/pr}dt)=C\frac1r(1-prt^{-1/pr}|_1^\infty)=C\frac{(1+pr)}r<\infty.
$$
For (iii), we note that if
 $\lim_{r\to\infty}\frac{1}{r}\tau(T^{p+\frac{1}{r}})$
exists, then $T^p\in {\mathcal Z}_1$ and by (i) $T\in{\mathcal
Z}_p$
\end{proof}

In view of the preceding corollary we have the following
implications
$$T\in {\mathcal L}^{p,\infty}\Longrightarrow T\in {\mathcal Z}_p,$$
$$T\in {\mathcal Z}_p\Longleftrightarrow T^p\in {\mathcal Z}_1={\mathcal L}^{1,\infty}.$$
Hence, everything which has been proved for $T\in {\mathcal
Z}_1={\mathcal L}^{1,\infty}$ is automatically true for $S=T^p$
provided $T\in {\mathcal Z}_p$ or especially if $T\in{\mathcal
L}^{p,\infty}$.
\subsection{The space ${\mathcal Z}_p, \ p>1$ is strictly larger than
${\mathcal L}^{p,\infty}$}
We deduce the result in the title of this subsection by proving
that the analogue of Theorem \ref{newresult} does not hold when
$p>1$.

\begin{prop} \label{newnegativeresult}
{\it The assumption $\sup_{r\ge
1}\frac{1}{r}\tau(T^{p+\frac{1}{r}})<\infty$ does not guarantee
$T\in\clL ^{(p,\infty)}$. }
\end{prop}
\begin{proof} We use the notation  $\mu_t(T):=x(t),\ t>0$.
The proof is based on the observation (see \cite{KPS} and also
detailed explanations in \cite[Section 5]{Suk}) that the ordinary
norm
$$
\|x\|_\psi=\sup_{t>0}\frac{\int_0^tx^*(s)ds}{\psi(t)}
$$
in the Marcinkiewicz space $M(\psi)$  (here, $\psi\in \Omega$ as
in Section 2) is equivalent to the quasi-norm
$$
F_\psi(x)=\sup_{0<t<\infty}\frac{tx^*(t)}{\psi(t)}
$$
provided that $\liminf_{t\to \infty}\frac{\psi(2t)}{\psi(t)}>1$.
For $\psi_p(t)=t^{1-1/p},\ p>1,$ the norm $\|\cdot\|_{\psi_p}$ and
quasi-norm $F_{p}(\cdot)=F_{\psi_p}(\cdot)$ are equivalent. In
other words, the norm of any element $T$ from the ideal $\clL
^{(p,\infty)}$ is equivalent to $F_p(x)$. This is not the case for
$\psi_0(t):=\ln(1+t)$ (that is the functional
$F_0(\cdot)=F_{\psi_0}(\cdot)$ and the norm in $\clL
^{(1,\infty)}$ are not equivalent) and it is easy to locate a
function $z(t)=z^*(t)$ such that $\|z\|_{\psi_0}<\infty$ but
$F_0(z)=\sup_{t>0}z^*(t)t=\infty$. For example, we take
$z(t)=n/2^{n^2}$ for $t\in (2^{(n-1)^2},2^{n^2}],\ n=1,2,...,$ and
$z(t)=1$ for $t\in [0,1]$. It is easy to verify that there exists
$0<C<\infty$ such that $$\int_0^tz^*(s)ds\le C\ln(1+t)$$ (that is
$z\prec\prec C/(1+t)$) and at the same time
$$z(t)t|_{t=2^{n^2}}=z(2^{n^2})2^{n^2}=n,\ n=1,2,...$$
(that is $F_0(z)=\infty$).

Observe that since $z\prec\prec C/(1+t)$, we have for every
$\nu>0$
$$
\int_0^\infty z(t)^{1+\nu}dt\le
C\int_0^\infty(1/(1+t))^{1+\nu}dt=C/\nu<\infty.
$$
Now, let us fix $p>1$ and set $x(t)=z^{1/p}(t)$ for $t>0$. The
estimate above gives
$$
s\int_0^\infty x^{p+s}(t)dt=p(s/p\int_0^\infty z(t)^{1+s/p} dt)\le
Cp<\infty.
$$
Nevertheless,
$$
F_p(x)=\sup_{0<t<\infty}{x(t)t^{1/p}}={(F_0(z))^{1/p}}=\infty.
$$
That is the condition $\sup_{r\ge
1}\frac{1}{r}\tau(T^{p+\frac{1}{r}})<\infty$ does not imply $T\in
\clL ^{(p,\infty)}$.
\end{proof}

We remark that while ${\mathcal Z}_p, \
p>1$ is the $p$-convexification of the ideal ${\mathcal
L}^{1,\infty}$; in turn, the ideal ${\mathcal L}^{p,\infty}$ is
the $p$-convexification of some subideal in ${\mathcal
L}^{1,\infty}$, which is termed the `small ideal' in \cite{CPS2}.
We will establish this latter
fact in subsection 5.2.
\subsection{Limits of zeta functions}
Our earlier results enable us to considerably weaken the hypotheses
in one of the main theorems of \cite{CPS2}.
First we recall the following preliminary result proved in \cite{CPS2}.
 \begin{prop}\label{2.2}(weak$^*$-Karamata theorem)
{\it Let $\tilde\omega\in L_\infty(\mbR)^*$ be a dilation
invariant state and let $\beta$ be a real valued, increasing,
right continuous function on $\mbR_+$ which is zero at zero and
such that the integral $h(r)=\int_0^\infty
e^{-\frac{t}{r}}d\beta(t)$ converges for all $r>0$ and
$C=\tilde\omega-\lim_{r\to\infty} \frac{1}{r}h(r)$ exists. Then}
$$\tilde\omega-\lim_{r\to\infty} \frac{1}{r}h(r)
=\tilde\omega-\lim_{t\to\infty}\frac{\beta(t)}{t}.$$
\end{prop}

The classical Karamata theorem has a similar statement with
the $\tilde\omega$ limits replaced by ordinary limits.

In the following we will take $T\in{\mathcal L}^{(1,\infty)}$ positive,
$||T||\leq 1$ with spectral resolution $T=\int \lambda dE(\lambda)$. We would
like to integrate with respect to $d\tau(E(\lambda))$; unfortunately, these
scalars $\tau(E(\lambda))$ are, in general,
 all infinite. To remedy this situation, we instead
must integrate with respect to the increasing (negative) real-valued function
$N_T(\lambda)=\tau(E(\lambda)-1)$ for $\lambda >0$. Away from $0$, the
increments $\tau(\triangle E(\lambda))$ and $\triangle N_T(\lambda)$ are, of
course, identical. The following theorem is a strengthened version of Theorem
3.1 of \cite{CPS2} made possible by Proposition \ref{prop2.4}.

\begin{thm}\label{newthm3.1}{\it
For $T\in{\mathcal L}^{(1,\infty)}$ positive, $||T||\leq 1$ %and $\tilde\omega\in L^\infty({\mathbb R})^*$
let $\omega$ be a $D_2$-dilation and $P^\alpha$-invariant,
$\alpha>1$ state on $L^\infty( \mathbb{R}^*_+)$. Let
$\tilde\omega=\omega\circ L$ where $L$ is given in Section 3, then
we have:
$$\tau_\omega(T)=\tilde\omega-\lim\frac{1}{r}\tau(T^{1+\frac{1}{r}}).$$
If $\lim_{r\to\infty}\frac{1}{r}\tau(T^{1+\frac{1}{r}})$ exists then
$$\tau_\omega(T)=\lim_{r\to\infty}\frac{1}{r}\tau(T^{1+\frac{1}{r}})$$
for an arbitrary dilation invariant functional
$\omega\in L^\infty({\mathbb R}^*_+)^*$.}
\end{thm}

\begin{proof} The proof is just a minor rewriting of the corresponding
argument in \cite{CPS2}. By  Proposition \ref{cps2prop}, the state
$\tilde\omega$ is dilation invariant and by Theorem
\ref{newresult}(i) $h(r)=\frac{1}{r}\tau(T^{1+\frac{1}{r}})\in
L^\infty(\mathbb{R}_+)$. So, we can apply the weak$^*$-Karamata
theorem. First write $\tau(T^{1+\frac{1}{r}})=\int_{0^+}^1
\lambda^{1+\frac{1}{r}} dN_T(\lambda)$. Thus setting
$\lambda=e^{-u}$
$$\tau(T^{1+\frac{1}{r}})=\int_0^\infty e^{-\frac{u}{r}}d\beta(u)$$
where $\beta(u)=\int_u^0 e^{-v}dN_T(e^{-v})=-\int_0^u e^{-v}dN_T(e^{-v})$.
Since the change of variable $\lambda=e^{-u}$ is strictly decreasing, $\beta$
is, in fact, nonnegative and increasing.
By the weak$^*$-Karamata theorem
 applied to $\tilde\omega\in L^\infty({\mathbb R})^*$
$$\tilde\omega-\lim_{r\to \infty} \frac{1}{r}\tau(T^{1+\frac{1}{r}})
=\tilde\omega-\lim_{u\to\infty}
\frac{\beta(u)}{u}.$$

Next with the substitution $\rho=e^{-v}$ we get:
\begin{equation}
\label{(3.2)}\tilde\omega-\lim_{u\to\infty}\frac{\beta(u)}{u}=
\tilde\omega-\lim_{u\to\infty}
\frac{1}{u}\int^1_{e^{-u}}
\rho dN_T(\rho).  \end{equation}
Set $f(u)=\frac{\beta(u)}{u}$.
We want to make the change of variable $u=\log t$
or in other words to consider $f\circ \log = Lf$.
This is permissable by the discussion in Section 3
which tells us that if we start with a %$G_2$ and $M$ invariant
functional $\omega\in L^\infty({\mathbb R}_+^*)^*$ as in the theorem
we may replace it by the functional $\tilde\omega=\omega\circ L$
which is dilation invariant
%is $G_1$ and $H$ invariant as required by the
%theorem.
with
$$\tilde\omega-\lim_{r\to \infty} \frac{1}{r}\tau(T^{1+\frac{1}{r}})
=\tilde\omega-\lim_{u\to\infty}\frac{\beta(u)}{u}$$
$$=\tilde\omega-\lim_{u\to\infty}
f(u)=\omega-\lim_{t\to\infty}Lf(t)=\omega-\lim_{t\to\infty}
\frac{1}{\log t}\int_{1/t}^1\lambda dN_T(\lambda).$$

Now, by Proposition \ref{prop2.4} and Corollary
\ref{after_prop2.4} applied to $\psi(t)=\log(1+t)\sim\log t$
$$\omega-\lim_{t\to\infty}
\frac{1}{\log t}\int_{1/t}^1\lambda dN_T(\lambda)
=\omega-\lim_{t\to\infty} \frac{1}{\log t}\tau(\chi_{(\frac{1}{t},1]}(T)T)
=\tau_\omega(T).$$
This completes the proof of the first part of the theorem.

The proof of the second part is similar. Using the classical
Karamata theorem
we obtain the following analogue of (\ref{(3.2)}):
$$\lim_{r\to\infty}\frac{1}{r}\tau(T^{1+r})=\lim\frac{\beta(u)}{u}=
\lim_{u\to\infty}\frac{1}{u}\int_{e^{-u}}^1\rho dN_T(\rho).$$
Making the substitution $u=\log t$ on the right hand side we have
 by Proposition \ref{prop2.4}
$$\lim_{u\to\infty}\frac{1}{u}\int_{e^{-u}}^1\rho dN_T(\rho)=
\lim_{t\to\infty}\frac{1}{\log t}\int_{\frac{1}{t}}^1\lambda
dN_T(\lambda)
=\tau_\omega(T)=\lim_{t\to\infty}\frac1{\log(1+t)}\int_0^t\mu_s(T)ds.$$
%where in the last equality we need only dilation invariance of the
%state $\omega\in L^\infty({\mathbb R}^*_+)^*$
\end{proof}

%Now we may use the discussion in this subsection and ideas from
%\cite{CPS2} to obtain some stronger consequences of a knowledge of
%the asymptotics of the zeta function. We use the discussion in
%Section~\ref{PrelimDilInvState} which tells us that if we start
%with an $M$ invariant functional $\omega\in L_\infty(\mbR_+^*)^*$
%then the functional $\tilde\omega=\omega\circ L$ is $H$ invariant.
%In addition we will assume $\omega$ is maximally invariant.

We now deduce some corollaries of the discussion above. Retaining
the notation as in the previous theorem we let $\omega$ be a
$D_2$-dilation and $P^\alpha$-invariant, $\alpha>1$ state on
$L^\infty( \mathbb{R}^*_+)$. Let $\tilde\omega=\omega\circ L$. The
assumption that $\frac{1}{r}\zeta(T^{1+\frac{1}{r}})$ is bounded
in $r$ means that, by Theorem \ref{newresult}, $T\in {\mathcal
Z}_1= {\mathcal L}^{1,\infty}$. Then by Theorem \ref{newthm3.1}
$$\tilde\omega-\lim_{r\to\infty} \frac{1}{r}\zeta(T^{1+\frac{1}{r}})
=\tau_\omega(T).$$
Consequently using (\ref{(+)}) if either $T\in
{\mathcal Z}_p$ or if $T\in {\mathcal L}^{p,\infty},\ p>1$ we have
the formulae
\begin{equation}\label{p>1}
\tilde\omega-\lim_{r\to\infty} \frac{1}{r}\zeta(T^{p+\frac{1}{r}})
=\tilde\omega-\lim_{r\to\infty}\frac1r\tau(T^{p+1/r})=p\tilde\omega-\lim_{pr\to\infty}\frac1{pr}\tau({T^p}^{(1+1/pr)})
=p\tau_\omega(T^p)
\end{equation}
where the last step uses dilation invariance of $\tilde\omega$,
which is guaranteed by our choice of $\omega$. The equation
(\ref{p>1}) together with Theorem \ref{newthm3.1} tell us that if
one of the limits in the previous equality is true then so are the
others. In particular, if $\lim_{r\to\infty}
\frac{1}{r}\zeta(T^{p+\frac{1}{r}})$ exists, then $T\in {\mathcal
Z}_p$ and
$$\lim_{r\to\infty} \frac{1}{r}\zeta(T^{p+\frac{1}{r}})=p\lim_{t\to\infty}\frac1{\log(1+t)}\int_0^t\mu_s(T^p)ds.$$
\section{The heat semigroup formula}\label{HSFsec}
\subsection{Asymptotics of the trace of the heat semigroup}
Throughout this section  $T\geq 0$. For $q\in {\mathbb R}_+$ we
define $e^{-T^{-q}}$ as the operator that is zero on $\ker T$ and
on $\ker T^\perp$ is defined in the usual way by the functional
calculus. We remark that if $T\geq 0$, $T\in{\mathcal Z}_{p}$ for
some $p\geq 1$ then  $e^{-tT^{-q}}$ is trace class for all $t>0$.
This is because if $x\in E$, where $(E, \|\cdot\|_E)$ is any
symmetric (or r.i.) space then
$$
\|x\|_E \ge \|x^*(t)\chi_{[0,s]}(t)\|_E \ge
x^*(s)\|\chi_{[0,s]}\|_E=x^*(s)\phi(s),
$$
where $\phi(\cdot)$ is the fundamental function of $E$.
Consequently, $x^*(s)\le \|x\|_E/\phi(s)$. For $E={\mathcal
Z}_p=\mathcal{L}_p^{1,\infty}$  (see Corollary \ref{OrM}(ii)) the
fundamental function is $\phi(s)=(s/log(1+s))^{1/p}$
%(see Corollary \ref{OrM} where it is shown that ${\mathcal Z}_p=\mathcal{L}_p^{1,\infty}$).
Hence, for every $t>0$
$$
\mu_s(e^{-tT^{-q}})=e^{-t/(\mu_s(T))^q}\le
e^{-tC(s/log(1+s))^{q/p}} \le e^{-tCs^{q/p-\epsilon}}$$ for some
$C>0$ all $0<p,q$ and $0<\epsilon<q/p$. Thus
$\tau(e^{-tT^{-q})}<\infty$ for $q>0$ (since $\epsilon>0$ is
arbitrary).

\begin{thm} (cf \cite{CPS2})\label{mainprop} {\it
If $T\geq 0$, $T\in{\mathcal Z}_{p}$, $1\leq p<\infty$ then,
choosing $\omega$ to be DPM invariant and $\tilde\omega$ to be
related with $\omega$
as in Remark 3.6, %{\bf (deleted $\epsilon$)}%and any $\epsilon>0$
we have
for $q>0$
$$\omega-\lim_{\lambda\to \infty} \frac{1}{\lambda}\tau(e^{-T^{-q}\lambda^{-q/p}})
=\frac1q\Gamma(p/q)\tilde\omega-\lim_{r\to\infty} \frac{1}{r}\zeta(p+\frac{1}{r})=\frac pq\Gamma(p/q)\tau_\omega(T^p).$$
}
\end{thm}

\begin{proof}
We have, using the Laplace transform,
$$T^s = \frac{1}{\Gamma(s/q)}\int_0^\infty t^{s/q -1}e^{-tT^{-q}}dt.$$
Then
$${\Gamma(s/q)}\zeta(s)={\Gamma(s/q)}\tau(T^s)=
\int_0^\infty t^{s/q -1}\tau(e^{-tT^{-q}})dt.$$
%Make the change of variable $t=1/\lambda^{2/p}$ so that the preceding
%formula becomes
%$$\frac{p}{2}\Gamma(s/2)\zeta(s)=
%\int_0^\infty \lambda^{-\frac{s}{p}-1}
%\tau(Ae^{-\lambda^{-2/p}T^{-2}})d\lambda.$$
We split this integral into two parts,
$\int_0^1$ and $\int_1^\infty$
and call the second integral $R(r)$ where $s=p+\frac{1}{r}$.
Then
$$R(r)=\int_1^\infty t^{p/q+1/(qr)-1}\tau(e^{-tT^{-q}})dt.$$
The integrand decays exponentially in $t$ as $t\to\infty$
because
$T^{-q}\geq {\Vert T^q\Vert}^{-1}\bf 1$
so that
$$\tau(e^{-tT^{-q}})\leq \tau(e^{-T^{-q}}e^{-\frac{t-1}{\Vert T^{q}\Vert}}).
$$
Then we can conclude that $R(r)$ is bounded independently of $r$ and so
$\lim_{r\to\infty}\frac{1}{r}R(r)=0$.
For the other integral
$\int_0^1 t^{p/q+1/(qr)-1}\tau(e^{-tT^{-q}})dt$ we can make the substitution
$t = e^{-\mu q/p}$. Then elementary calculus gives
$$\int_0^1 t^{p/q+1/(qr)-1}\tau(e^{-tT^{-q}})dt=
-q/p\int_\infty^0 e^{-\mu(1+\frac1{pr}))}\tau(e^{-e^{-\mu q/p} T^{-q}})d\mu
=q/p\int_0^\infty e^{-\frac{\mu}{pr}} d\beta(\mu)$$
where
$\beta(\mu)=\int_0^\mu e^{-v}\tau(e^{-e^{-vq/p}T^{-q}})dv$.
Hence we can now write
$$\Gamma(p/q+\frac{1}{rq})\zeta(p+\frac{1}{r})
=q/p\int_0^\infty e^{-\frac{\mu}{pr}} d\beta(\mu)+R(r).$$

Then we have (remembering that the term $\frac{1}{r}R(r)$ has limit
zero as $r\to \infty$)
%Since $\lim_{u\to 1}\Gamma(u)=1$ we have
$$\tilde\omega-\lim_{r\to\infty}\frac{1}{r}
\Gamma(p/q+\frac1{pr})\zeta(p+\frac{1}{r})=
\Gamma(p/q)\tilde\omega-\lim_{r\to\infty}\frac{1}{r}\zeta(p+\frac{1}{r})$$
$$=\tilde\omega-\lim_{r\to\infty}\frac{q}{pr}
\int_0^\infty e^{-\mu/pr}d\beta(\mu)
=q\tilde\omega-\lim_{r\to\infty}\frac{1}{r}
\int_0^\infty e^{-\mu/r}d\beta(\mu)$$
where the last step uses the
assumed dilation invariance of $\tilde\omega$.
So
$$\tilde\omega-\lim_{r\to\infty}\frac{1}{r}
\Gamma(p/q+\frac1{pr})\zeta(p+\frac{1}{r})=
q \tilde\omega-\lim_{r\to\infty}\frac{1}{r}\int_0^\infty e^{-\frac{\mu}{r}} d\beta(\mu)
$$

Now we are exactly in a position to use the weak*-Karamata
theorem above to evaluate the RHS. Indeed, we now conclude
$$\tilde\omega-\lim_{r\to\infty}\frac{1}{r}
\int_0^\infty e^{-\frac{\mu}{r}} d\beta(\mu)
=\tilde\omega-\lim_{\mu\to\infty}
 \frac{\beta(\mu)}{\mu}.$$
We can summarise the preceding in the equation
\begin{equation}\label{(++)}
\Gamma(p/q)\tilde\omega-\lim_{r\to\infty}\frac{1}{r}
\zeta(p+\frac{1}{r})=q\tilde\omega-\lim_{\mu\to\infty}
 \frac{\beta(\mu)}{\mu}.
\end{equation}
Now make the change of variable $\lambda = e^v$
in the defining expression for $\beta(\mu)$
to obtain
$$\frac{\beta(\mu)}{\mu}=\frac{1}{\mu}\int_1^{e{^\mu}}
\lambda^{-2}\tau(e^{-\lambda^{-q/p}T^{-q}})d\lambda$$

Make the substitution $\mu =\log t$ so the RHS becomes
$$\frac{1}{\log
t}\int_1^t\lambda^{-2}\tau(e^{-T^{-q}\lambda^{-q/p}})d\lambda=g_1(t)$$
This is the Cesaro mean of
$$g_2(\lambda)=\frac{1}{\lambda}\tau(e^{-T^{-q}\lambda^{-q/p}}).$$
Thus
 as we chose $\omega\in L_\infty( \mathbb{R}^*_+)^*$ to be $M$ invariant
and $\tilde\omega$ to be related to $\omega$ as in Remark
\ref{pairs_of_functionals}
we have
$$\tilde\omega-\lim_{\mu\to\infty}
 \frac{\beta(\mu)}{\mu}=\omega(g_1)=\omega(g_2).$$
%
%$$\Gamma(p/q)\tilde\omega-\lim_{r\to\infty}
%\frac{1}{r}\zeta(p+\frac{1}{r}) =p\Gamma(p/q)\tau_\omega(T^p)$$
Then using (\ref{(++)}), we obtain
$$\Gamma(p/q)\tilde\omega-\lim_{r\to\infty}
\frac{1}{r}\zeta(p+\frac{1}{r})=
q\omega(g_2)=q\omega-\lim_{\lambda\to \infty}
\frac{1}{\lambda}\tau(e^{-T^{-q}\lambda^{-q/p}})$$ Thus by
(\ref{p>1}) we obtain the statement of the theorem:
$$\Gamma(p/q)\tilde\omega-\lim_{r\to\infty}
\frac{1}{r}\zeta(p+\frac{1}{r})=
q\omega-\lim_{\lambda\to \infty} \frac{1}{\lambda}\tau(e^{-T^{-q}\lambda^{-q/p}})
=p\Gamma(p/q)\tau_\omega(T^p).$$
\end{proof}

%To prove the theorem take $A\geq 0$ and use
%the Proposition \ref{mainprop} to assert that
%$$\frac{p}{q}\Gamma(p/q)\tau_\omega(AT^p)
%%\Gamma()\tilde\omega-\lim_{r\to\infty}\frac{1}{r}\tau(AT^{1+\frac{1}{r}})
%=\omega-\lim_{\lambda\to\infty}\lambda^{-1}
%\tau(Ae^{-\lambda^{-q/p}T^{-q}}).
%$$
%Then for self adjoint $A$
%write $A=a^+-a^-$ where $a^\pm$ are positive so that
%$$\Gamma(1+p/q)\tau_\omega(AT^p)=
%\Gamma(1+p/q)(\tau_\omega(a^+T)-\tau_\omega(a^-T))$$
%$$=
%\omega-\lim_{\lambda\to\infty}\lambda^{-1}
%\tau(a^+e^{-\lambda^{-q/p}T^{-q}})-
%\omega-\lim_{\lambda\to\infty}\lambda^{-1}\tau(a^-e^{-\lambda^{-q/p}T^{-q}})
%$$
%$$=
%\omega-\lim_{\lambda\to\infty}\lambda^{-1}\tau(Ae^{-\lambda^{-q/p}T^{-q}}).$$
%We can extend to general bounded $A$ by a similar argument.
\subsection {The $L^{p,\infty}$-case and the 'small' ideal.}
As $T\in {\mathcal L}^{p,\infty}$ means that
$\mu_t(T)t^{1/p}<C<\infty$ and $\mu_t(T^p)=\mu_t(T)^p$ we conclude that
$
\mu_t(T^p)t<C^p<\infty.$
That is
$T\in {\mathcal L}^{p,\infty}\Longrightarrow S=T^p\in {\mathcal I}$ where
 ${\mathcal I}$ is the so called `small' subideal
of ${\mathcal L}^{1,\infty}$ identified in \cite{CPS2}.
Recall that $\mathcal I$
is specified by the condition on the singular values
of $T\geq 0, T\in{\mathcal L}^{1,\infty} $: $\mu_s(T)\leq C/s$ for some
constant $C>0$.
In subsection 4.1 \cite{CPS2} we proved the following result by
a direct argument that avoids the use of the zeta function.
If $\omega$ is $M$ invariant and satisfies conditions (1),(2),(3) of Theorem
3.4 and $T\in \mathcal I$ then
$$\omega-\lim_{\lambda\to\infty}
\lambda^{-1}\tau(e^{-\lambda^{-2}T^{-2}})=\Gamma(3/2)\tau_\omega(T).$$

We may now apply this stronger result of \cite{CPS2}
to operators $S\in \mathcal I$ where $S=T^p$
and  $T\in {\mathcal L}^{p,\infty}$
to obtain the equality
$$\omega-\lim_{\lambda\to\infty}
\lambda^{-1}\tau(e^{-\lambda^{-2}S^{-2}})=\Gamma(3/2)\tau_\omega(S).$$
Hence we obtain the following result
\begin{equation}\label{(!)}
\mbox{\it If $T\in {\mathcal L}^{p,\infty}$ then }
\omega-\lim_{\lambda\to\infty} \lambda^{-1}\tau(e^{-\lambda^{-2}T^{-2p}})=\Gamma(3/2)\tau_\omega(T^p).
\end{equation}}
Note that we have obtained this result under weaker conditions on
$\omega$ than the more general Theorem \ref{mainprop} where
$T\in{\mathcal Z}_p$. It would be interesting to understand an
example in noncommutative geometry where ${\mathcal Z}_p$ arises
naturally. We remark that in classical geometric examples
such as differential operators on manifolds it is $
{\mathcal L}^{p,\infty}$ $p\geq 1$ and the `small ideal' $\mathcal I$
that arise naturally.

  A further idea motivated by the geometric case is that one may argue
the other way, from a knowledge of the asymptotics of the trace of
the heat semigroup, to information on the zeta function.
Thus let us assume that the
trace of the heat operator $\tau(e^{-tT^{-2}})$ exists for all
$t>0$ and in addition
has an asymptotic expansion in inverse powers of $t$ as $t\to 0$.
These assumptions hold for Dirac Laplacians for example
in classical geometry
and it is well known
in this case that one can infer from the asymptotic expansion
the nature of the first singularity of
$\zeta(s)$ (as $\Re s$ decreases)
 from the leading term in inverse powers of $t$.
We now explain this in some detail.

Thus assume that $\tau(e^{-tT^{-2}})= Ct^{-p/2}$ + lower order
powers of $t^{-1}$ as $t\to 0$. We recall that as in  Theorem \ref{mainprop}
$\tau(e^{-tT^{-2}})\to 0$ exponentially as $t\to\infty$.
We introduce
\begin{equation}\label{zeta}
\zeta_1(s)=
\frac{1}{\Gamma(s/2)}\int_0^1 t^{s/2 -1}\tau(e^{-tT^{-2}})dt,\  s>p
\end{equation}
and $$\zeta_2(s)= \frac{1}{\Gamma(s/2)}\int_1^\infty t^{s/2
-1}\tau(e^{-tT^{-2}})dt, s>0.$$ Then $\zeta_2$ is analytic in a
neighborhood of $s=p$ and
we may write $\zeta(s)=\tau(T^s):=\zeta_1(s)+
\zeta_2(s)$ for $\Re s>p$. Then the only contribution to the
singularity at $s=p$ comes from $\zeta_1$. Now
$$\frac{1}{\Gamma(s/2)}\int_0^1 t^{s/2 -1} Ct^{-p/2} dt=
\frac{C}{\Gamma(s/2)(s/2-p/2)}$$ and thus substitution in (\ref{zeta})
gives
$$\zeta(s)= \tau(T^s)= \frac{C}{\Gamma(s/2)(s/2-p/2)} + K(s)$$
where $K(s)$ is holomorphic for $s=p$.
(We note that the lower order terms in the asymptotic expansion
do contribute to the term $K(s)$ but these contributions
are analytic near $s=p$.)
Thus we may take the limit $\lim_{s\to p}(s-p)\zeta(s)$ and only
the first term contributes as $\lim_{s\to p}(s-p)K(s)=0$.

\begin{prop}
If  $\tau(e^{-tT^{-2}})$ has an asymptotic expansion in inverse
powers of $t$ with the leading term being $C/t^{p/2}$ for some
constant $C$ then $T\in{\mathcal Z}_{p}$ and $$\lim_{s\to
p}(s-p)\tau(T^s)=p\tau_\omega(T^p)$$ for any
 $D_2$ (and $M$) invariant $\omega$.
\end{prop}

\section{Application to spectral triples}

Throughout this Section the following assumptions hold. We let
$\mathcal D$ be an unbounded self adjoint densely defined operator
on $\mathcal H$ affiliated to $\mathcal N$ (this amounts to
$(1+\mathcal D^2)^{-1}\in \mathcal N$). We suppose that $\mathcal
A$ is a *-algebra in $\mathcal N$ consisting of operators $a$ such
that $[\clD ,a]$ is bounded and  refer to the triple
$(\clD,\clA,{\mathcal N})$ as a semifinite spectral triple.

Denote for brevity $\mathcal M^\psi:=M(\psi)({\mathcal N},\tau)$
with $\psi$ as in Section 4, satisfying (\ref{psicon}). As in
Corollary \ref{OrM}, we consider the following $p$-convexification
of $\mathcal M^\psi$
$$
{\mathcal M}^{\psi,p}:=\{T\in {\mathcal N}_+:\
\|T\|_{\psi,p}=\sup_{1<u<\infty}\frac{(\int_0^u
\mu_t(T)^pdt)^{1/p}}{\psi^{1/p}(u)}<\infty\},\quad p>1.
$$
We let $\tau_\omega$ be a Dixmier trace on $\mathcal M^\psi$
corresponding to a suitable singular state $\omega$. Suppose that
$(1+ \mathcal D^2)^{-p/2}\in \mathcal M^\psi$, or equivalently
that $(1+ \mathcal D^2)^{-1/2}\in \mathcal M^{\psi,p}$. In
applications of noncommutative geometry the functional
$\phi_\omega$ on $\mathcal A$ given by $\phi_\omega(a)=
\tau_\omega(a(1+D^2)^{-p/2})$ plays a key role. In particular it
is of interest to know if this functional is a trace on $\mathcal
A$. In \cite{CGS} this question was answered in the affirmative
for the case of $(1+ \mathcal D^2)^{-1/2}\in {\mathcal
L}^{p,\infty} $. Their proof generalizes to our setting. In
particular, it holds under the weaker assumption $(1+ \mathcal
D^2)^{-1/2}\in\mathcal Z_p$.

\begin{thm} {\it Under the immediately preceding hypotheses
we have} $$\phi_\omega(ab)= \phi_\omega(ba)  \ \ \ \ \ a,b \in
\mathcal A.$$
\end{thm}

The proof is an extension of the approach in \cite{CGS}. We need
four preliminary facts. Some may be proved in a similar way to the
corresponding results in \cite{CGS}.

\begin{lemma} {\it  Given a spectral triple $(\clD,\clA,{\mathcal N})$ we have\\
(i) For $a,b\in \mathcal N$ the H\"older inequality
$$\tau_\omega(ab)\leq\tau_\omega(|a|^p)^{1/p}\tau(|b|^q)^{1/q}$$
for $p,q\geq 1$, $\frac{1}{p}+\frac{1}{q} = 1$, holds.\\
(ii) For any $r$ with $0<r<1$ and $a \in \mathcal A$
the operator $[(1+\clD^2)^{r/2},a]$ is bounded and satisfies
$$||[(1+\clD^2)^{r/2},a]||\leq C ||[\clD,a]||$$
where the constant $C>0$ does not depend on $a$.\\
(iii) Let $T\in \mathcal M^\psi$ and  $f(t)=\mu_t(T)$ so
that $f$ is a bounded decreasing function on $(0,\infty)$ from
$M(\psi)$, then $f^\alpha \in L_1({\mathbb R}_+)$ for every
$\alpha >1$.\\
(iv) The statement of the theorem (for
$(1+\clD^2)^{-p/2}\in \mathcal M^\psi$) is implied by  }
$$\tau_\omega(|[(1+\clD^2)^{-p/2},a]|)=0 \mbox{ for all } a \in {\mathcal A}.$$
\end{lemma}

\begin{proof} (i)  We have by \cite[Proposition 1.1]{CS} and by the
H\"older inequality for function spaces
$$\int_0^t \mu_s(ab) ds \leq \int_0^t \mu_s(a)\mu_s(b) ds \leq
(\int_0^t  \mu_s(a)^p ds)^{1/p}(\int_0^t  \mu_s(b)^q ds)^{1/q}$$
Dividing by $\psi(t)$ and applying
the functional $\omega$ we get
$$\tau_\omega(ab)\leq \omega\left [\left(\frac{\int_0^t  \mu_s(a)^p ds}{\psi(t)}\right)^{1/p}
\left(\frac{\int_0^t  \mu_s(b)^q
ds}{\psi(t)}\right)^{1/q}\right]$$
$$\leq  \omega\left(\frac{\int_0^t  \mu_s(a)^p ds}{\psi(t)}\right)^{1/p}
\omega \left(\frac{\int_0^t  \mu_s(b)^q ds}{\psi(t)}\right)^{1/q}=
\tau_\omega(|a|^p)^{1/p}\tau_\omega(|b|^q)^{1/q}$$ using H\"older
inequality for states on abelian $C^*$-algebras. We omit the proof
for $p=1$, $q=\infty$.

(ii) If $\mathcal N$ is taken in its left regular representation,
then the claim follows immediately from \cite[Theorem 3.1]{PS}.
The general case is done in \cite[Theorem 2.4.3]{P}. Note, that
the assumption  made in \cite{CGS} that $\mathcal D$ has a bounded
inverse is now redundant.

(iii) Using the inequalities preceding Lemma 4.3, we have for any
$\beta>1$ $f(t) \leq C'\frac{1}{t ^{1/\beta}}$ for some $C'>0$ and
all sufficiently large $t$'s. Since $\alpha>1$ is given, we can
choose $\beta$ so that $\frac{\alpha}{\beta} =\gamma >1$, and so
$f^\alpha(t) \leq
C'/{t ^{\gamma}}$ which gives the required result.

(iv) Let $T=(1+\clD^2)^{-p/2}$ and $a,b\in \mathcal A$. Then we
know that for $T'\in \mathcal M^\psi$
 $\tau_\omega(T'a)=\tau_\omega(aT')$ (see \cite{Co4} or
\cite[Lemma 3.2(i)]{CPS2}) and hence
$$\phi_\omega([a,b])=\tau_\omega(Tab-aTb)=\tau_\omega([T,a]b).$$
Then
$$|\tau_\omega([T,a]b)|\leq \tau_\omega(|[T,a]|)||b||)=0$$
with the last equality is implied by the hypothesis of the lemma.
\end{proof}

Choose $r$ with $0<r<1$ such that $k=p/r \in \mathbb N$.
Following \cite{CGS}, we see that the proof of the theorem rests
on the identity (for $k\in\mathbb N$)
$$[a,(1+\clD^2)^{-kr/2}]=\sum_{j=1}^k (1+\clD^2)^{-jr/2}[(1+\clD^2)^{r/2},a]
(1+\clD^2)^{(j-k-1)r/2}$$ where we are using part (ii) of the
Lemma to give boundedness of $[(1+\clD^2)^{r/2},a]$. We now apply
the previous identity to obtain:
$$\tau_\omega(|[a,(1+\clD^2)^{-p/2}]|)= \tau_\omega(|[a,(1+\clD^2)^{-kr/2}]|)$$
$$
\leq \sum_{j=1}^k \tau_\omega[ |(1+\clD^2)^{-jr/2}[(1+\clD^2)^{r/2},a]
(1+\clD^2)^{(j-k-1)r/2}|]$$
Hence choosing $p_j= \frac{2p}{r(2j-1)}$, $q_j= \frac{2p}{r(2k-2j+1)}$
and applying part (i) of the Lemma,
$$\tau_\omega(|[a,(1+\clD^2)^{-p}]|)\leq ||[(1+\clD^2)^{r/2},a]||
 \sum_{j=1}^k (\tau_\omega((1+\clD^2)^{-p_jjr/2}))^{1/p_j}
(\tau_\omega((1+\clD^2)^{(j-k-1)q_jr/2}))^{1/q_j}$$
The exponents $p_jjr/2$ and $(j-k-1)q_jr/2$ are larger than $p$ so
using part (iii) of the Lemma, the Dixmier trace
in the last two terms vanishes.
Now use part (iv) of the Lemma to complete the proof of the Theorem.

\def\itakdalee{$\dots$}

\end{document}